\documentclass[a4paper,reqno]{amsart}
\usepackage{pdfsync}
\usepackage{amsmath,amssymb}
\usepackage{enumerate}
\usepackage{xspace}
\usepackage{boxedminipage}
\usepackage{listings}% for including source code
\usepackage{tabularx}
\usepackage{subfigure}
\usepackage{algpseudocode}
\usepackage{graphicx}
       \newtheoremstyle%
           {plain}% name
           {}% Space before, vuoto = `valore di default'
           {}% Space after
           {\mdseries\slshape}% body font
           {}% Indent (empty = no indent, \parindent = para indent)
           {\bfseries}% Thm head font
           {.}% Punctuation after the heading
           {1.0ex}% Space after heading: \newline = to start at next line
           {}% Thm head spec (can be left empty, meaning `normal')
           
       \newtheoremstyle
           {note}% name
           {}% Space before, vuoto = `valore di default'
           {}% Space after
           {}% body font
           {}% Indent (empty = no indent, \parindent = para indent)
           {\bfseries}% Thm head font
           {.}% Punctuation after the heading
           {1.0ex}% Space after heading: \newline = to start at next line
           {}% Thm head spec (can be left empty, meaning `normal')
           
       \swapnumbers
       \theoremstyle{plain}
       \theoremstyle{note}
       \newenvironment{Example}[1][]{\subsection{Example\xspace{\ifx&#1&{}\else{ (#1)}\fi}}}{}
       \newenvironment{Rem}[1][]{\subsection{Remark\xspace{\ifx&#1&{}\else{ (#1)}\fi}}}{}
       \newenvironment{Obs}[1][]{\subsection{Remark\xspace{\ifx&#1&{}\else{ (#1)}\fi}}}{}
       \newenvironment{Defn}[1][]{\subsection{Definition\xspace{\ifx&#1&{}\else{ (#1)}\fi}}}{}
       
       \newenvironment{Hyp}[1][]{\subsection{Assumption\xspace{\ifx&#1&{}\else{ (#1)}\fi}}}{}

              \providecommand{\Le}[1]{\ensuremath{\leb{#1}}} 
              \providecommand{\hoz}{\sobhz1(\W)}

              \providecommand{\laginterpol}[1]{\La^{#1}}%{\mathcal I^{#1}}

              \providecommand{\tol}{\ensuremath{\operatorname{tol}}\xspace}

              \providecommand{\T}[1]{\cT^{#1}}

              \providecommand{\ie}{i.e.,\xspace}

              \providecommand{\tr}[1]{\ensuremath{\operatorname{tr}{#1}}}

              \providecommand{\rangefromto}[3]{\ensuremath{#1=#2,\ldots,#3}}
              \providecommand{\figwidth}{\textwidth} \providecommand{\figscale}{1}

              \providecommand{\fes}{\ensuremath{\mathbb V}}
              \renewcommand{\sin}[1]{\ensuremath{\operatorname{sin}\left(#1\right)}}
              \renewcommand{\cos}[1]{\ensuremath{\operatorname{cos}\left(#1\right)}}

              \providecommand{\MA}{Monge--Amp\`ere }
              \providecommand{\cof}[1]{\operatorname{Cof} #1}
              \renewcommand{\det}[1]{\operatorname{det} #1}

              \providecommand{\frob}[2]{\ensuremath{{#1}{:}{#2}}}
              \renewcommand{\div}[1]{\operatorname{div}\left(#1\right)}
              \renewcommand{\H}{\ensuremath{\vec{H}}}
              
              \providecommand{\symm}{\Symmatrices d}
              \providecommand{\feszero}{\mathring{\fes}}
              \providecommand{\Nzero}{\mathring{N}}

              \providecommand{\A}{\ensuremath{\vec{A}}}
              \providecommand{\N}{\ensuremath{\vec{N}}}

              \providecommand{\matentry}[3]{{#1}_{#2,#3}}
              \renewcommand{\matentry}[3]{{#1}_{#2,#3}}
              \providecommand{\nlo}{\cN}
              \providecommand{\nlop}[1]{\nlo[#1]}
              \providecommand{\nlfunk}{F}
              \providecommand{\Dnlop}[1]{\D\cN[#1]}
              \providecommand{\MAD}{Monge--Amp\`ere--Dirichlet\xspace}

              \newcommand{\eig}[1]{\sigma\qp{#1}}
        \usepackage{ifthen}
        \usepackage{amsmath}
        \usepackage{amssymb}    %mathematics
        \usepackage{stmaryrd} %more symbols
        \usepackage{eufrak}     %fraktur charcters
        \usepackage{mathrsfs} %nice calligraphic characters
        \usepackage{bbold}    %blackboard bold
        \provideboolean{usemathrsfs}
        \setboolean{usemathrsfs}{true}
        \usepackage{enumerate}% a facility to customize numbered lists
        \usepackage{xspace}% smart spacing after macros: detects punctuation
        \usepackage[greek,british,english]{babel}
        \usepackage{verbatim}% provides also the {comment} environment
        \usepackage{fancyvrb}
        \usepackage{hyperref}
        \usepackage{listings}% for colorized code output
        \usepackage{xcolor}
        \provideboolean{usebeamer}% when using beamertex, 
        \provideboolean{isthesis}% used in the PhD thesis
        \provideboolean{isamsltex}% use when using AMS papers style
        \provideboolean{issiamltex}% on when using SIAM papers... becoming obsolete
    \usepackage{xcolor}
    \colorlet{a}{red}
    \colorlet{b}{blue}
    \colorlet{c}{green!50!blue}
    \colorlet{d}{magenta}
    \colorlet{e}{cyan}
    \colorlet{f}{yellow!50!black}
    \colorlet{g}{white}
    \colorlet{h}{black!50}
    \colorlet{i}{black}
    \colorlet{j}{black!75} 
    \providecommand{\linkedurl}[1]{\url{#1}}
    \providecommand{\linkedemail}[1]{\href{mailto:#1}{#1}}
    \providecommand{\email}[1]{{\linkedemail{#1}}}
    \providecommand{\Ignore}[1]{}
    \providecommand{\ignore}[1]{}
    \providecommand{\freeze}[1]{}% Obsolete, kept for compatibility
    \providecommand{\crossout}[1]{{\textcolor{red!20}{#1}}}
    \providecommand{\highlight}[1]{{\color{blue}#1}}
    \providecommand{\standout}[1]{\colorbox{a}{\textcolor{g}{#1}}}
    
    \provideboolean{shownotes}
    \setboolean{shownotes}{true}%%false
    \newcounter{margnote}[page]
    \providecommand{\margnotemark}{{\standout{\upshape\texttt{\arabic{margnote}}}}}
    \providecommand{\margnote}[2][]{
      \ifthenelse{
        \boolean{shownotes}
      }{
        \stepcounter{margnote}
        \margnotemark
        \marginpar{
          \texttt{
            \begin{minipage}{2cm}
              \raggedright\tiny
              \margnotemark
              {\ifx&#1&{}\else{#1 says:}\fi} 
              #2
            \end{minipage}
          }
        }
      }{
      }
    }
    \providecommand{\mathnote}[2][]{
      \ifthenelse{
        \boolean{shownotes}
      }{
        \stepcounter{margnote}
        \text{
          \standout{
            \texttt{
              \tiny
              \margnotemark
              #1:
              #2
            }
          }
        }
      }{
      }
    }
    \provideboolean{showtodo}
    \providecommand{\todo}[1]{\ifthenelse{\boolean{showtodo}}{\margnote[To do.]{#1}}{}}
    \providecommand{\Todo}[1]{
      \ifthenelse{\boolean{showtodo}}{
        \begin{center}
        \begin{tikzpicture}
         \node[fill=a!17]{
           \begin{minipage}{\textwidth}
             \texttt{\bfseries{\small #1}}
           \end{minipage}
         };
        \end{tikzpicture}
        \end{center}
      }{}}

    \provideboolean{showcomments}
    \providecommand{\margincomment}[1]{
    \ifthenelse{\boolean{showcomments}}{\marginpar{\tiny #1}}{}
    }
    \provideboolean{showchanges}
    \setboolean{showchanges}{false}
    \providecommand{\changes}[1]{
      \ifthenelse{\boolean{showchanges}}{{\highlight{#1}}}{#1}
    }
    \providecommand{\changefromto}[3][replace with]{
      \ifthenelse{\boolean{showchanges}}
      {{\crossout{#2}\margnote{#1}}{\highlight{#3}}}
      {#3\xspace}
    }
    \providecommand{\ChangePar}[2]{
      \ifthenelse{\boolean{showchanges}}
      {{\par$\mapsfrom$ \textcolor{red!20}{#1}}{\par$\mapsto$ \textcolor{blue}{#2}}}
      {\par #2}
    }
    \providecommand{\InsertPar}[1]{
      \ifthenelse{\boolean{showchanges}}
      {{\par$\mapsto$ \textcolor{blue}{#1}}}
      {\par #1}
    }

    \usepackage{mathrsfs}%extra mathematical symbols
    \provideboolean{usemathrsfs}
    \setboolean{usemathrsfs}{true}% switch if mathrsfs is unused
    \providecommand{\mathscript}
    	   {\mathscr}
     \providecommand{\cA}{\ensuremath{\mathscript A}\xspace}
     \providecommand{\cB}{\ensuremath{\mathscript B}\xspace}
     \providecommand{\cC}{\ensuremath{\mathscript C}\xspace}

     \providecommand{\cM}{\ensuremath{\mathscript M}\xspace}
     \providecommand{\cN}{\ensuremath{\mathscript N}\xspace}

     \providecommand{\cT}{\ensuremath{\mathscript T}\xspace}

     \providecommand{\bbbold}{\mathbb}

     \providecommand{\rN}{\ensuremath{\bbbold N}\xspace}
     
     \providecommand{\rP}{\ensuremath{\bbbold P}\xspace}
     
     \providecommand{\rR}{\ensuremath{\bbbold R}\xspace}

     \providecommand{\rV}{\ensuremath{\bbbold V}\xspace}

    \providecommand{\ie}{\ensuremath{\text{ i.e., }}\xspace}
    
     \providecommand{\naturals}{\rN\xspace}
     
     \providecommand{\NO}{\ensuremath{\rN_0}}

     \providecommand{\reals}{\rR}
     \providecommand{\R}[1]{\reals^{#1}}

     \providecommand{\closure}[1]{\overline{#1}}

     \providecommand{\La}{\ensuremath{\varLambda}\xspace}

     \providecommand{\W}{\ensuremath{\varOmega}\xspace}
     \providecommand{\qp}[1]{\ensuremath{\!\left({#1}\right)}}
     \providecommand{\qpreg}[1]{\ensuremath{(#1)}}
     \providecommand{\qpbig}[1]{\ensuremath{\big(\!#1\!\big)}}
     \providecommand{\qpBig}[1]{\ensuremath{\Big(\!#1\!\Big)}}
     \providecommand{\qpbigg}[1]{\ensuremath{\bigg(\!\!#1\!\!\bigg)}}
     \providecommand{\qpBigg}[1]{\ensuremath{\Bigg(\!\!#1\!\!\Bigg)}}
     \providecommand{\qb}[1]{\ensuremath{\!\left[{#1}\right]}}

     \providecommand{\qa}[1]{\ensuremath{\left\langle{#1}\right\rangle}}

     \providecommand{\powqp}[2]{\ensuremath{\qp{#2}^{\kern -.25em\lower .2ex\hbox{\scriptsize $#1$}}\kern-.1em}}
     \providecommand{\powqpreg}[2]{\ensuremath{\qpreg{#2}^{\kern -.2em\lower .3ex\hbox{\scriptsize $#1$}}\kern-.3em}}
     \providecommand{\powqpbig}[2]{\ensuremath{\qpbig{#2}^{\kern -.2em\lower .3ex\hbox{\scriptsize $#1$}}\kern-.3em}}
     \providecommand{\powqpBig}[2]{\ensuremath{\qpBig{#2}^{\kern -.2em\lower .3ex\hbox{\scriptsize $#1$}}\kern-.3em}}
     \providecommand{\powqpbigg}[2]{\ensuremath{\qpbigg{#2}^{\kern -.2em\lower .3ex\hbox{\scriptsize $#1$}}\kern-.3em}}
     \providecommand{\powqpBigg}[2]{\ensuremath{\qpBigg{#2}^{\kern -.2em\lower .3ex\hbox{\scriptsize $#1$}}\kern-.3em}}

     \providecommand{\norm}[1]{\ensuremath{\left|#1\right|}}
     
     \providecommand{\abs}[1]{\ensuremath{\left|#1\right|}}
     \providecommand{\Norm}[1]{\ensuremath{\left\|#1\right\|}}

     \providecommand{\ltwop}[2]{\ensuremath{\qa{#1,#2}}}

     \providecommand{\duality}[3][]{\ensuremath{#1\langle #2\,#1\vert\,#3#1\rangle}}
     \providecommand{\ensemble}[2]{\ensuremath{\left\{ #1:\;#2 \right\}}}

     \providecommand{\squmatmd}[3][1]{
       \begin{bmatrix}
         \matentry{#2}{#1}{#1}&\dotsc&\matentry{#2}{#1}{#3}
         \\
         \vdots & \ddots &\vdots
         \\
         \matentry{#2}{#3}{#1}&\dotsc&\matentry{#2}{#3}{#3}
       \end{bmatrix}
     }
     
     \providecommand{\sumifromto}[3]{\ensuremath{\sum_{#1=#2}^{#3}}}

     \providecommand{\matentry}[3]{\bientry{\mat{#1}}{#2}{#3}}

     \providecommand{\rangefromto}[3]{\ensuremath{#1\integerbetween{#2}{#3}}}
    
     \providecommand{\D}{\ensuremath{\,\mathrm{D}}}
    
    \providecommand{\registered}%
    {\ensuremath{^\text{\textregistered}}}
    
    \providecommand{\AND}{\ensuremath{\text{ and }}}

    \providecommand{\OR}{\ensuremath{\text{ or }}}

    \providecommand{\constant}[1]{\ensuremath{C_{#1}}}
    \providecommand{\constref}[2][]{\ensuremath{\constant{\text{\ref{#2}{\ifx&#1&{}\else{,#1}\fi}}}}}
    \renewcommand{\div}{\operatorname{div}}

    \providecommand{\eye}{\operatorname{\bf I}}               % identity matrix
    \providecommand{\Eye}[1]{
      \begin{bmatrix}
      \ifthenelse{#1>1}{
        \ifthenelse{#1>2}{
          \ifthenelse{#1>3}{
            1&0&\dotso&0
            \\
            0&1&\dotso&0
            \\
            \vdots&\vdots&\ddots&\vdots
            \\
            0&0&\dotso&1
          }{
            1&0&0
            \\
            0&1&0
            \\
            0&0&1
          }
        }{
          1&0
          \\
          0&1
        }
      }{
        1
      }
      \end{bmatrix}
    }
    \providecommand{\Oh} {\operatorname{O}}                   % Landau's Big O
    \providecommand{\pd}[1]{\ensuremath{\partial_{#1}}\xspace} % basic
    \providecommand{\pdt}[1][]{\pd t{{\ifx&#1&{}\else{\qb{#1}}\fi}}}                       % Partial derivative in t
    \providecommand{\trace}{\operatorname{trace}}             % trace (of a matrix)

    \providecommand{\transpose}{{\boldsymbol\intercal}}   % transpose symbol
    
    \providecommand{\Transpose}[1]{\ensuremath{{#1}^{\transpose}}}
    
    \providecommand{\Transposevec}[1]{\Transpose{\vec{#1}}}
    \providecommand{\transposevec}[1]{\Transposevec{#1}}

    \providecommand{\Hess}{\ensuremath{\D^2}}
    \providecommand{\intersected}{\ensuremath{\cap}}
    \providecommand{\meet}{\intersected}

    \providecommand{\union}[1]{\ensuremath{\bigcup}_{#1}}

    \renewcommand{\vec}[1]{\ensuremath{\boldsymbol{#1}}}
    \providecommand{\geovec}[1]{\vec{#1}}
    \providecommand{\geomat}[1]{\vec{#1}}

    \providecommand{\boundary}{\partial}
     \providecommand{\CC}{\ensuremath{\operatorname C}\xspace}%Continuous functions
     \providecommand{\HH}{\ensuremath{\operatorname H}\xspace}
     \providecommand{\LL}{\ensuremath{\operatorname L}\xspace}
     
     \providecommand{\WW}{\ensuremath{\operatorname W}\xspace}

     \providecommand{\cont}[1]{\ensuremath{\CC^{#1}}}

     \providecommand{\leb}[1]{\ensuremath{\LL_{#1}}}
     \providecommand{\lebloc}[1]{\ensuremath{{\LL^{\mathrm{loc}}_{#1}}}}
    
     \providecommand{\sob}[2]{\ensuremath{{\smash\WW}^{#1}_{#2}}}

     \providecommand{\sobh}[1]{\ensuremath{\HH^{#1}}}
     \providecommand{\sobhz}[1]{\sobh{#1}_0}
     
     \providecommand{\poly}[1]{\ensuremath{\rP}^{#1}}
    
    \providecommand{\Symmatrices}[1]{\ensuremath{\operatorname{Sym}{(\R{d\times d})}}}

     \providecommand{\fespace}{\rV}
     \providecommand{\fezerospace}{\ensuremath{\smash{\mathring\fespace}}}
     
     \providecommand{\fes}[1]{\ensuremath{\fespace^{#1}}}
     \providecommand{\fez}[1]{\ensuremath{\fezerospace^{#1}}}

    \providecommand{\Forall}{\:\forall\:}
    
    \providecommand{\Foreach}{\quad\Forall}

    \usepackage{amscd}

    \providecommand{\funk}[3]{\ensuremath{#1:#2\to#3}}

    \providecommand{\restriction}[2]{\left.#1\right|_{#2}}
    \renewcommand{\restriction}[2]{\left.#1\right|_{#2}}

    \providecommand{\Algoname}[1]{\ensuremath{\text{\textsf{#1}\xspace}}}

    \providecommand{\euro}{\textgreek{\euro}}
    
    \providecommand{\ListParameters}{}
    \renewcommand{\ListParameters}% these can be changed by renewcommand
    {
    	 \setlength{\topsep}{0em}
    	 \setlength{\leftmargin}{0em}
             \setlength{\itemsep}{0ex}
    	 \setlength{\parsep}{.5ex}
    	 \setlength{\itemindent}{\labelsep}
    	 \addtolength{\itemindent}{\labelwidth}
    }
    
    \newcounter{LetterListItem}
    \renewcommand{\theLetterListItem}{(\alph{LetterListItem})}
    {
    	\begin{list}%
    	{\theLetterListItem\ }%
    	{\usecounter{LetterListItem}
    	 \ListParameters
    	}
    }%
    {\end{list}}
    \newcounter{NumberListItem}
    \renewcommand{\theNumberListItem}{\arabic{NumberListItem}}
    {
    	\begin{list}%
    	{\theNumberListItem.\ }%
    	{\usecounter{NumberListItem}%
    	 \ListParameters
    	}
    }%
    {\end{list}}
    \newcounter{QuestionListItem}
    \renewcommand{\theQuestionListItem}{\textbf{Question \arabic{QuestionListItem}}}
    {
    	\begin{list}%
    	{\theQuestionListItem.\ }%
    	{\usecounter{QuestionListItem}%
    	 \ListParameters
    	}
    }%
    {\end{list}}
    \newcounter{RomanListItem}
    \renewcommand{\theRomanListItem}{(\roman{RomanListItem})}
    {
    	\begin{list}%
    	{\theRomanListItem\ }%
    	{\usecounter{RomanListItem}
    	 \ListParameters
    	}
    }%
    {\end{list}}
    \newcounter{StepsItem}
    {
    	\begin{list}%
    	{Step \theStepsItem.\ }%
    	{\usecounter{StepsItem}%
    	 \ListParameters
    	}
    }%
    {\end{list}}
    
    \providecommand{\grad}{\nabla}
    \renewcommand{\grad}{\nabla}
    
    \usepackage{amsthm} % an extra package for theorem-like environments
    \ifthenelse{\boolean{isthesis}}{%
      \setcounter{secnumdepth}{1}%
    }{
      \setcounter{secnumdepth}{2} % this number is the depth of numbering
    }
    \providecommand{\ListParameters}{}
    \renewcommand{\ListParameters}
    {
    	 \setlength{\topsep}{0em}
    	 \setlength{\leftmargin}{0em}
             \setlength{\itemsep}{0ex}
    	 \setlength{\parsep}{.5ex}
    	 \setlength{\itemindent}{\labelsep}
    	 \addtolength{\itemindent}{\labelwidth}
    }
    \newtheoremstyle{plain}% name
      {}% Space before, vuoto = `valore di default'
      {}% Space after
      {\mdseries\slshape}% body font
      {\parindent}% Indent (empty = no indent, \parindent = para indent)
      {\bfseries}% Thm head font
      {.}% Punctuation after the heading
      {.5em}% Space after heading: \newline = to start at next line
      {}% Thm head spec (can be left empty, meaning `normal')
    
    \newtheoremstyle{note}% name
      {}% Space before
      {}% Space after
      {}% body font
      {\parindent}% Indent (empty = no indent, \parindent = para indent)
      {\bfseries}% Thm head font
      {.}% Punctuation after the heading
      {.5em}% Space after heading: \newline = to start at next line
      {}% Thm head spec (can be left empty, meaning `normal')
    
    \newtheoremstyle{claim}% name
      {}% Space before, vuoto = `valore di default'
      {}% Space after
      {\mdseries\slshape}% body font
      {}% Indent (empty = no indent, \parindent = para indent)
      {\bfseries}% Thm head font
      {}% Punctuation after the heading
      {.5em}% Space after heading: \newline = to start at next line
      {}% Thm head spec (can be left empty, meaning `normal')
    
    \newtheoremstyle{exercise}% name
      {}%  Space before, empty = `default'
      {}% Space after
      {}% body font
      {}% Indent (empty = no indent, \parindent = para indent)
      {\bfseries}% Thm head font
      {.}% Punctuation after the heading
      {1em}% Space after heading: \newline = to start at next line
      {}% Thm head spec (can be left empty, meaning `normal')
    
    \newtheoremstyle{break}% name
      {}%  Space before, empty = `default'
      {}% Space after
      {}% body font
      {}% Indent (empty = no indent, \parindent = para indent)
      {\bfseries}% Thm head font
      {.}% Punctuation after the heading
      {\newline}% Space after heading: \newline = to start at next line
      {}% Thm head spec (can be left empty, meaning `normal')
    
      \providecommand{\SolName}{Solution}%Solution       %Soluzione    %Lösung
      \providecommand{\Proofname}{Proof}%Dimostrazione   %Démonstration%Beweis
    
    \ifthenelse{\boolean{isthesis}}{%
      
    }{%
      
    }
    \providecommand{\pdfformat}[1]{
       \provideboolean{pdfoutput}
       \setboolean{pdfoutput}{#1}% default is dvi and ps output
      \ifthenelse{\boolean{pdfoutput}}{
        \typeout{using pdf}
        \providecommand{\graphext}{pdf}
        \renewcommand{\graphext}{pdf}
        \providecommand{\graphextex}{pdf_t}
        \renewcommand{\graphextex}{pdf_t}
      }{
        \typeout{using eps}
        \usepackage[dvips]{graphicx,xcolor}
        \providecommand{\graphext}{eps}
        \renewcommand{\graphext}{eps}
        \providecommand{\graphextex}{eps_t}
        \renewcommand{\graphextex}{eps_t}
      }
      \usepackage{epsfig}
      \usepackage{tikz}
      \usepackage{rotating}
      \definecolor{SussexFlint}{rgb}{.00,.19,.21}
      \definecolor{SussexGrey}{rgb}{.51,.58,.49}
      \definecolor{SussexOrange}{rgb}{.94,.29,.00}
      \definecolor{SussexYellow}{rgb}{1.00,.73,.00}
      \definecolor{SussexRed}{rgb}{.94,.01,.49}
      \definecolor{SussexPurple}{rgb}{.48,.06,.44}
      \definecolor{SussexGreen}{rgb}{.00,.58,.46}
      \definecolor{SussexBlue}{rgb}{.00,.58,.65}
      \colorlet{a}{SussexOrange}
      \colorlet{b}{SussexRed}
      \colorlet{c}{SussexYellow}
      \colorlet{d}{SussexPurple}
      \colorlet{e}{SussexGreen}
      \colorlet{f}{SussexBlue}
      \colorlet{g}{white}%%backGround
      \colorlet{h}{SussexGrey}%%haze
      \colorlet{i}{black}%%ink
      \colorlet{j}{SussexFlint}
      \newcommand{\mausDarkColorTheme}{
        \colorlet{a}{SussexYellow!50!yellow}
        \colorlet{b}{SussexGreen!50!green}
        \colorlet{c}{SussexBlue}%%{SussexBlue!50!cyan}
        \colorlet{d}{SussexOrange!50!yellow}
        \colorlet{e}{SussexRed!50!red}
        \colorlet{f}{SussexPurple!50!magenta}
        \colorlet{g}{black}%% backGround
        \colorlet{h}{SussexFlint!50!black}
        \colorlet{i}{white}%% Ink color
        \colorlet{j}{SussexGrey}
      }
    }
    \providecommand{\solution}{\textbf{\SolName.}\xspace}
    
     \newcounter{phantombox}[enumi]%%[equation]
     \provideboolean{showphantoms}
     \renewcommand{\thephantombox}{\roman{phantombox}}
     \newcommand{\phantombox}[1]{\stepcounter{phantombox}
       \ensuremath{\boxed{
           \ifthenelse
    	   {\boolean{showphantoms}}
    	   {#1^{\phantom{\textup{(\thephantombox)}}}}
    	   {\phantom{#1}^{\textup{(\thephantombox)}}}
         }
       }
     }

     \provideboolean{hidesolution}
     \newcommand{\consolution}[1]{
       \ifthenelse{
         \boolean{hidesolution}
       }{
       }{
         {\par \small {\solution}\ #1\par\ \\[5pt]}}
     }
     \provideboolean{showmarks}
     \renewcommand{\marks}[1]{
       \ifthenelse{\boolean{showmarks}}{\marginpar{{\tiny [$#1$ marks]}}}{}}
     \newcommand{\condibreak}{\ifthenelse{\boolean{hidesolution}}{\newpage}{}}
     
\usepackage{natbib}
\usepackage[utf8]{inputenc}
\usepackage{tikz}
\providecommand{\nlop}[1]{\cN[#1]}
\providecommand{\nlfunk}{F}
\providecommand{\Dnlop}{\mat N}

\numberwithin{equation}{section}
\setboolean{showtodo}{true}%false
\setboolean{shownotes}{true}%false
\setboolean{showchanges}{true}%false
\setboolean{usemathrsfs}{false}%true
\providecommand{\fenics}{FEniCS\xspace}

\providecommand{\paraview}{ParaView\xspace}
\newcommand{\gnuplot}{Gnuplot\xspace}
\newcommand{\dolfin}{DOLFIN\xspace}

\author{Omar Lakkis}
\address{ Omar Lakkis\newline
  Department of Mathematics\newline
  University of Sussex\newline
  Brighton, England\newline
  GB-BN1 9RF} \curraddr{}
\email{\linkedemail{o.lakkis@sussex.ac.uk}}
\urladdr{\linkedurl{http://www.maths.sussex.ac.uk/Staff/OL}}
\author{Tristan Pryer} 
\address{Tristan Pryer\newline
  School of Mathematics, Statistics and Actuarial Sciences\newline 
  University of Kent\newline
  Canterbury\newline 
  UK-CT2 7NF, United Kingdom} 
\curraddr{}
\email{\linkedemail{t.pryer@kent.ac.uk}}
\title[FEM for nonlinear elliptic problems]{A finite element method
  for nonlinear elliptic problems}
\date{\today}

\begin{document}

\maketitle
\begin{abstract}
  We present a continuous finite element method for some examples of
  fully nonlinear elliptic equation.  A key tool is the discretisation
  proposed in Lakkis \& Pryer (2011) allowing us to work directly on
  the strong form of a linear PDE. An added benefit to making use of
  this discretisation method is that a \emph{recovered (finite
    element) Hessian} is a biproduct of the solution process. We build
  on the linear basis and ultimately construct two different
  methodologies for the solution of second order fully nonlinear PDEs.
  Benchmark numerical results illustrate the convergence properties of
  the scheme for some test problems as well as the Monge--Ampère
  equation and the Pucci equation.
\end{abstract}

%%%%%%%%%%%%%%%%%%%%%%%%%%%%%%%%%%%%%%%%%%%%%%%%%%%%%%%%%%%%%%%%%%%%%%%%
      %%% tex-main-file: "../nonlinear.tex" ***
      %%%š
      \section{Introduction}
      \label{sec:intro}
      
      Fully nonlinear PDEs arise in many areas, including differential
      geometry (the \MA equation), mass transportation (the
      Monge--Kantorovich problem), dynamic programming (the Bellman
      equation) and fluid dynamics (the geostrophic equations).  The
      computer approximation of the solutions of such equations is thus an
      important scientific task.  There are at least three main difficulties
      apparent to someone attempting to derive numerical methods for fully
      nonlinear equations: first, the strong nonlinearity on the highest
      order derivative which generally precludes a variational formulation,
      second, a fully nonlinear equation does not always admit a classical
      solution, even if the problem data is smooth, and the solution has to
      sought in a generalised sense (e.g., viscosity solutions), which is
      bound to slow down convergence rates, and third, a common problem in
      nonlinear solvers, the exact solution may not be unique and
      constraints, such as convexity requirements must be included in the
      constraints to ensure uniqueness.
      
      Regardless of the problems, the \emph{numerical approximation of fully
        nonlinear second order elliptic equations}, as described in
      \cite{CaffarelliCabre:1995}, have been the object of considerable recent
      research, particularly for the case of \MA of which
      \cite{OlikerPrussner:1988, LoeperRapetti:2005, DeanGlowinski:2006,
        FengNeilan:09:vanishing, Oberman:2008, Awanou:10, DavydovSaeed:12:techreport,
        BrennerGudiNeilanSung:11, Froese:11} are selected examples.
      
      For more general classes of fully nonlinear equations some methods
      have been presented, most notably, at least from a theoretical view
      point, in \cite{Bohmer:2008} where the author presents a $\cont{1}$
      finite element method shows stability and consistency (hence
      convergence) of the scheme, following a classical ``finite
      difference'' approach outlined by \cite{Stetter:1973} which requires a
      high degree of smoothness on the exact solution. From a practical
      point of view this approach presents difficulties, in that the
      $\cont{1}$ finite elements are hard to design and complicated to
      implement, in \cite{DavydovSaeed:10} a useful overview of
      Bézier-Bernestein splines in two spatial dimensions is provided and a
      full implementation in \cite{DavydovSaeed:12:techreport}. Similar
      difficulties are encountered in finite difference methods and the
      concept of \emph{wide-stencil} appears to be useful, for example by
      \cite{KuoTrudinger:92,KuoTrudinger:05,Oberman:2008,Froese:11}.
      
      In \cite{FengNeilan:09:mixed,FengNeilan:09:vanishing,Awanou:10}
      the authors give a method in which they approximate the general second
      order fully nonlinear PDE by a sequence of fourth order quasilinear
      PDEs. These are quasilinear biharmonic equations which are discretised
      via mixed finite elements, or using high-regularity elements such as
      splines.  In fact for the \MA equation, which admits two solutions, of
      which one is convex and another concave, this method allows for the
      approximation of both solutions via the correct choice of a
      parameter. On the other hand although computationally less expensive
      than $\cont{1}$ finite elements (an alternative to mixed methods for
      solving the biharmonic problem), the mixed formulation still results
      in an extremely large algebraic system and the lack of maximum
      principle for general fourth order equations makes it hard to apply
      vanishing viscosity arguments to prove convergence.  A somewhat
      different approach, based on $C^0$-penalty, has been recently proposed
      by \cite{BrennerGudiNeilanSung:11}, as well as ``pseudo time'' one by
      \cite{Awanou:11}.
      
      It is worth citing also a \emph{least square} approach described by
      \cite{DeanGlowinski:2006}.  This method consists in minimising the
      mean-square of the residual, using a Lagrange multiplier method. Also
      here a fourth order elliptic term appears in the energy.
      
      In this paper, we depart from the above proposed methods and explore a
      more ``direct'' approach by applying the \emph{nonvariational finite
        element method}, introduced in \cite{LakkisPryer:2011}, as a solver
      for the Newton iteration directly derived from the PDE.  To be more
      specific, consider the following model problem
      \begin{equation}
          \label{eq:modelproblem}
          \nlop{u} := \nlfunk(\Hess u) - f = 0
      \end{equation}
      with homogeneous Dirichlet boundary conditions where $\funk f\W\reals$
      is prescribed function and $\funk F\symm\reals$ is a real-valued
      algebraic function of symmetric matrixes, which provides an elliptic
      operator in the sense of \cite{CaffarelliCabre:1995}, as explained
      below in Definition \ref{def:ellipticity}.  The method we propose,
      consists in applying a Newton's method, given below by equation
      \eqref{eq:linearisedequation} of the PDE \eqref{eq:modelproblem},
      which results in a sequence of linear nonvariational elliptic PDEs
      that fall the framework of the nonvariational finite element method
      (NVFEM) proposed in \cite{LakkisPryer:2011}.  The results in this
      paper are computational, so despite not having a complete proof of
      convergence, we test our algorithm various problems that are
      specifically constructed to be well posed.  In particular, we test our
      method on the \MA problem, which is the de-facto benchmark for
      numerical methods of fully nonlinear elliptic equations.  This is in
      spite of \MA having an extra complication, which is conditional
      ellipiticity (the operator is elliptic only if the function is convex
      or concave.  A crucial, empirically observed feature of our method is
      that the convexity (or concavity) is automatically preserved if one
      uses $\poly2$ elements or higher.  For $\poly1$ elements this is not
      true and the scheme must be stabilized by reenforcing convexity (or
      concavity) at each timestep.  This was achieved in \cite{Pryer:2010}
      using a semidefinite programming method.  In a different spirit, but
      somewhat reminiscent, a stabilization procedure was obtained in
      \cite{BrennerGudiNeilanSung:11} by adding a penalty term.
      
      The rest of this paper is set out as follows. In
      \S\ref{sec:notation_and_discretisation} we introduce some notation,
      the model problem, discuss its ellipticity and Newton's method, which
      yields a sequences of nonvariational linearised PDE's. In
      \S\ref{sec:ndfem} we review of the nonvariational finite element
      method proposed in \cite{LakkisPryer:2011} and apply it to discretise
      the nonvariational linearised PDE's in Newton's method. In
      \S\ref{sec:unconstrained-fnl-pde} we numerically demonstrate the
      performance of our discretisation on a class of fully nonlinear PDE,
      those that are elliptic and well posed without constraining our
      solution to a certain class of functions.  In \S\ref{sec:ma} we turn
      to conditionally elliptic problems by dealing with the prime example
      of such problems, i.e., \MA.  We apply the discretisation to the \MA
      equation making use of the work \cite{Aguilera:2008} to check
      \emph{finite element convexity} is preserved at each iteration.
      Finally in \S\ref{sec:Pucci} we address the approximation of Pucci's
      equation, which is another important example of fully nonlinear
      elliptic equation.
      
      All the numerical experiments for this research, were carried out
      using the \dolfin interface for \fenics \cite{LoggWells:2010} and
      making use of \gnuplot and \paraview for the graphics.
      %%%%%%%%%%%%%%%%%%%%%%%%%%%%%%%%%%%%%%%%%%%%%%%%%%%%%%%%%%%%%%%%%%%%%%%%
      %%%
      %%% Local Variables: ***
      %%% mode:latex ***
      %%% tex-main-file: "../nonlinear.tex"  ***
      %%% End: ***
%%%%%%%%%%%%%%%%%%%%%%%%%%%%%%%%%%%%%%%%%%%%%%%%%%%%%%%%%%%%%%%%%%%%%%%
      %%%%%%%%%%%%%%%%%%%%%%%%%%%%%%%%%%%%%%%%%%%%%%%%%%%%%%%%%%%%%%%%%%%%%%%%
      %%% mode:latex ***
      %%% tex-main-file: "../nonlinear.tex"  ***
      %%%%%%%%%%%%%%%%%%%%%%%%%%%%%%%%%%%%%%%%%%%%%%%%%%%%%%%%%%%%%%%%%%%%%%%%
      \section{Notation}
      \label{sec:notation_and_discretisation}
      \subsection{Functional set-up}
      Let $\W\subset \R{d}$ be an open and bounded Lipschitz domain. We
      denote $\leb{2}(\W)$ to be the space of square (Lebesgue) integrable
      functions on $\W$ together with {its} inner product
      $\ltwop{v}{w} := \int_\W v w$ and norm $\Norm{v} :=
      \Norm{v}_{\leb2(\W)} = \ltwop{v}{v}^{1/2}$. %\margnote{Introduce $\cD$} 
      We denote by
      $\duality{v}{w}$ the action of a distribution $v$ on the function
      $w$. %If both $v,w\in\Le{2}(\W)$ then $\duality{v}{w} =
      %\ltwop{v}{w}$. We also denote by $\langle f \rangle_\w$ the integral
      %of a function $f$ over the domain $\w$ and drop the subscript for $\w
      %= \W$.
      
      We use the convention that the derivative $\D u$ of a function
      $u:\W\to\reals$ is a row vector, while the gradient of $u$, $\nabla u$
      is the derivatives transpose (an element of $\reals^d$, representing
      $\D u$ in the canonical basis). Hence
      \begin{equation}
        \nabla u 
        =
        \Transpose{\left(\D u\right)}.
      \end{equation}
      For second derivatives, we follow the common innocuous abuse of
      notation whereby the Hessian of $u$ is denoted as $\Hess u$ (instead
      of the more consistent $\D\grad  u$) and is represented by a $d\times d$
      matrix.
      
      The standard Sobolev spaces are ~\cite{Ciarlet:1978,Evans:1998}
      \begin{gather}
        \sobh{k}(\W)
        := 
        \sob{k}{2}(\W) 
        = 
        \ensemble{\phi\in\leb2(\W)}
                 {\sum_{\norm{\vec\alpha} \leq k}
                   \D^{\vec\alpha}\phi\in\leb2(\W)},
                 \\
        \hoz := \text{closure of }\cont{\infty}_0(\W) \text{ in } \sobh{1}(\W)
      \end{gather}
      where $\vec\alpha = \{ \alpha_1,...,\alpha_d\}$ is a
      multi-index, $\norm{\vec\alpha} = \sum_{i=1}^d\alpha_i$ and
      derivatives $\D^{\vec\alpha}$ are understood in a weak sense. 
      
      We consider the case when the model problem \eqref{eq:modelproblem} is
      uniformly elliptic in the following sense.
      
      \begin{Defn}[{ellipticity \cite{CaffarelliCabre:1995}}]
        \label{def:ellipticity}
        The operator $\nlop\cdot$ in Problem (\ref{eq:modelproblem}) is 
        called \emph{elliptic on $\cC\subseteq\symm$} if and only if
        for each $\geomat M\in\cC$ there exist $\Lambda\geq\lambda>0$, that
        may depend on $\geomat M$ such that
        \begin{equation}
          \label{eqn:conditional-ellipticity-of-nl}
          \lambda \sup_{\norm{\geomat \xi} = 1}\norm{\geomat N\geomat \xi}
          \leq 
          \nlfunk(\geomat M+\geomat N) - \nlfunk(\geomat M) 
          \leq
          \La \sup_{\norm{\geomat \xi} = 1}\norm{\geomat N\geomat \xi}
          \Foreach \geomat N\in\symm.
        \end{equation}
        If the largest possible set $\cC$ for which
        (\ref{eqn:conditional-ellipticity-of-nl}) is satisfied is a proper
        subset of $\symm$ we say that $\nlop\cdot$ is \emph{conditionally elliptic}.
        
        The operator $\nlop\cdot$ in Problem (\ref{eq:modelproblem}) is called
        to be \emph{uniformly elliptic} if and only if for some
        $\lambda,\La>0$, called \emph{ellipticity constants}, we have
        \begin{equation}
          \label{eq:ellipticity-of-nl}
          \lambda \sup_{\norm{\geomat \xi} = 1}\norm{\geomat N\geomat \xi}
          \leq 
          \nlfunk(\geomat M+\geomat N) - \nlfunk(\geomat M) 
          \leq
          \La \sup_{\norm{\geomat \xi} = 1}\norm{\geomat N\geomat \xi}
          \Foreach \geomat N,\geomat M\in\symm.
        \end{equation}
      \end{Defn}
      %%%%%%%%%%%%%%%%%%%%%%%%%%%%%%%%%%%%%%%%%%%%%%%%%%%%%%%%%%%%%%%%%%%%%%%%
      %%\subsection{Smooth elliptic}
      
      If $\nlfunk$ is differentiable (\ref{eq:ellipticity-of-nl}) can be
      obtained from conditions on the derivative of $\nlfunk$. A generic
      $\geomat M\in\symm$ is written as
      \begin{equation}
        \geomat M =
        \squmatmd md
        ;
      \end{equation}
      so the derivative of $\nlfunk$ in the direction $\geomat N$ is given by
      \begin{equation}
        \D\nlfunk(\geomat M)\geomat N
        =
        \frob{\nlfunk'(\geomat M)}{\geomat N}
      \end{equation}
      where the \emph{derivative matrix} $\nlfunk'(\geomat M)$ is defined by
      \begin{equation}
        \nlfunk'(\geomat M)
        :=
        \squmatmd{\partial\nlfunk(\geomat M)/\partial m}{d}
        %%  the JOY of MACROS
        %%  \begin{bmatrix}
        %%  \fracl{\partial F(\geomat M)}{\partial m_{1,1}} 
        %%  &
        %%  \dots
        %%  & \fracl{\partial F(\geomat M)}{\partial m_{d,1}} 
        %%  \\
        %%  \vdots & \ddots  & \vdots\\
        %%  \fracl{\partial F(\geomat R)}{\partial r_{1,d}}
        %%  & \dots  &
        %%  \fracl{\partial F(\geomat R)}{\partial r_{d,d}} 
        %%  \end{bmatrix}
        .
      \end{equation}
      %%%%%%%%%%%%%%%%%%%%%%%%%%%%%%%%%%%%%%%%%%%%%%%%%%%%%%%%%%%%%%%%%%%%%%%%
      %%\begin{Pro}[ellipticity criterion for differentiable elliptic
      %%    operators] 
      
            %% The ellipticity condition
            %% \begin{equation}
            %%   \tag{EC}
            %%   \label{eqn:ellipticity-of-nl}
            %% \end{equation}
            %% can be simplified for \emph{smooth nonlinearities}.
            %%%%%%%%%%%%%%%%%%%%%%%%%%%%%%%%%%%%%%%%%%%%%%%%%%%%%%%%%%%%%%%%%%%%%%%%
            Suppose $\nlfunk$ is differentiable. Then
            (\ref{eqn:conditional-ellipticity-of-nl}) is satisfied if and only if
            for each $\geomat M\in\cC$ there exists $\mu>0$ such that
            \begin{equation} 
              \label{eq:ellipticity-for-derivative}
              \Transpose{\geovec\xi}\nlfunk'(\geomat M)\geovec\xi 
              \geq 
              \mu
              \norm{\geovec \xi}^2 
              \Foreach 
              %% \geomat M\in\symm,
              \geovec\xi\in\R d.  
            \end{equation}
            Furthermore $\cC=\symm$ and $\mu$ is independent of $\geomat M$
            if and only if (\ref{eq:ellipticity-of-nl}) is satisfied.
      %%\end{Pro}
      %%%%%%%%%%%%%%%%%%%%%%%%%%%%%%%%%%%%%%%%%%%%%%%%%%%%%%%%%%%%%%%%%%%%%%%%
      \begin{Hyp}[smooth elliptic operator]
        \label{hyp:smooth-elliptic-operator}
        In the remainder of this paper we shall assume that
        $\nlop\cdot$ is conditionally elliptic on $\cC$ and
        \begin{equation}
          \nlfunk \in \cont{1}(\cC).
        \end{equation}
        Unless otherwise stated we will also assume that $\cC=\symm$.
      \end{Hyp}
      \subsection{Newton's method}
      \label{sec:newtons-method}
      The smoothness assumption \ref{hyp:smooth-elliptic-operator} allows
      to apply Newton's method to solve Problem (\ref{eq:modelproblem}).
      
      Given the initial guess $u^0\in\cont2(\W)$, with $\Hess u^0\in\cC$, for
      each $n\in\NO$, find $u^{n+1}\in\cont2(\W)$ with $\Hess u^{n+1}\in\cC$ such that
      \begin{equation}
        \label{eq:newtonsmethod}
        \Dnlop{u^n}\left(u^{n+1} - u^n\right) = -\nlop{u^n},
      \end{equation}
      where $\Dnlop u$ indicates the (Fréchet) derivative, which is formally
      given by
      \begin{equation}
        \label{eq:newtonsmethod2}
        \begin{split}
          {\Dnlop u}{v}
          &=
          \lim_{\epsilon\rightarrow 0} 
          \frac{\nlop{u + \epsilon v} - \nlop{u}}{\epsilon}
          \\ 
          &= 
          \lim_{\epsilon\rightarrow 0}
          \frac{\nlfunk(\Hess u + \epsilon \Hess v) - \nlfunk(\Hess u)}{\epsilon}
          \\ 
          &=
          \nlfunk'(\Hess u) : \Hess v,
        \end{split}
      \end{equation}
      for each $v\in\cont2(\W)$.
      Combining (\ref{eq:newtonsmethod}) and (\ref{eq:newtonsmethod2}) then
      results in the following nonvariational sequence of linear PDEs. Given
      $u^0$ for each $n\in\naturals_0$ find $u^{n+1}$ such that
      \begin{equation}
        \label{eq:system-of-nonvar-pdes}
        \nlfunk'(\Hess u^n) : \Hess\left(u^{n+1}-u^n\right) = f - \nlfunk(\Hess u^n).
      \end{equation}
      %An abstract summary of Newton's method on Banach spaces can be found
      %in \S\ref{sec:newtons-method-on-banach-spaces}.
      
      The PDE (\ref{eq:system-of-nonvar-pdes}) comes naturally in a nonvariational
      form. If we attempted to rewrite into a variational form, in order, say, to
      apply a ``standard'' Galerkin method, we would introduce an advection term
      which would depend on derivatives of $\nlfunk'$, \ie for generic $v,w$
      \begin{equation}
        \frob{\nlfunk'(\Hess v)}{\Hess w}
        =
        \div\qb{\nlfunk'(\Hess v) \nabla w}- \div\qb{ \nlfunk'(\Hess v) }\nabla w.
      \end{equation}
      where the matrix-divergence is taken row-wise:
      \begin{equation}
        \div\qb{\nlfunk'(\Hess v(\vec x))} 
        :=
        \qp{
          \sum_{i=1}^d \frac{\partial}{\partial x_i} 
          \qb{\matentry{F'}i1(\Hess v(\vec x))}
          , 
          \dotsc
          ,
          \sum_{i=1}^d \frac{\partial}{\partial x_i} 
          \qb{\matentry{F'}id(\Hess v(\vec x))}
          %%    \sum_{i=1}^d \frac{\partial}{\partial x_i} \qp{\matentry{\partial\nlfunk(\geomat M)/\partial m}{i}{d}}
        }
      \end{equation}
      and the chain rule provides us, for each $\rangefromto j1d$, with
      \begin{equation}
          \sum_{i=1}^d \frac{\partial}{\partial x_i} 
          \qb{\matentry{F'}ij(\Hess v(\vec x))}
          =
          \sumifromto{k,l}1d
          \partial_{k,l}\matentry{F'}ij
          (\Hess v(\vec x))
          \partial_{ikl}v(\vec x).
      \end{equation}
      This procedure is undesirable for many reasons. Firstly it requires
      $F$ to be twice differentiable and it involves a third order
      derivative of the functions $u^{n+1}$ and $u^n$ appearing in
      (\ref{eq:newtonsmethod}).  Moreover, the ``variational'' reformulation
      could very well result in the problem becoming advection dominated and
      unstable for conforming FEM, as was manifested in numerical examples
      for the linear equation \cite[\S 4.2]{LakkisPryer:2011}.  In order to
      avoid these problems, we here propose the use of the nonvariational
      finite element method described next.
      %%%%%%%%%%%%%%%%%%%%%%%%%%%%%%%%%%%%%%%%%%%%%%%%%%%%%%%%%%%%%%%%%%%%%%%%
      %%%%%%%%%%%%%%%%%%%%%%%%%%%%%%%%%%%%%%%%%%%%%%%%%%%%%%%%%%%%%%%%%%%%%%%%
      \section[The NVFEM]{The nonvariational finite element method}
      \label{sec:ndfem}
      
      %% The structure of 
      %% (\ref{eq:system-of-nonvar-pdes}) motivates the use of the
      %% nonvariational finite element method (NVFEM) introduced in
      %% \cite{LakkisPryer:2011}.
      In \cite{LakkisPryer:2011} we proposed the \emph{nonvariational finite
        element method} (NVFEM) to approximate the solution of problems of
      the form (\ref{eq:system-of-nonvar-pdes}).  We review here the NVFEM
      and explain how to use it in combination with the Newton method to
      derive a practical Galerkin method for the numerical approximation of
      Problem (\ref{eq:modelproblem})'s solution.
      %%%%%%%%%%%%%%%%%%%%%%%%%%%%%%%%%%%%%%%%%%%%%%%%%%%%%%%%%%%%%%%%%%%%%%%%
      \subsection{Distributional form of (\ref{eq:system-of-nonvar-pdes}) and generalised Hessian}
      Let $\A \in \Le{\infty}(\W)^{d\times d}$ and for each $\geovec
      x\in\W$, let $\A(\geovec x)\in\symm$, the space of bounded, symmetric,
      positive definite, $d \times d$ matrices and $\funk f\W\reals$. The
      \emph{Dirichlet linear nonvariational elliptic problem} associated
      with $\A$ and $f$ is
      \begin{equation}
        \frob{\A}\Hess u=f
        \AND
        \restriction u{\boundary\W}=0.
      \end{equation}
      Testing this equation, and assuming $u\in\sobh2(\W)\cap\hoz$ such
      that $\restriction{\grad u}{\boundary\W}\in\leb2(\boundary\W)$,
      we may write it as
      \begin{equation}
        \begin{split}
          \label{Problem}
          \ltwop{\frob{\A}{\Hess u}}{\phi} 
          &= \ltwop{f}{\phi}
          \qquad \Foreach
          \phi\in\cont\infty_0(\W).
        \end{split}
      \end{equation}
      To allow a Galerkin type discretisation of (\ref{Problem}), we need to
      restrict the test functions $\phi$ to finite element function spaces
      that are generally \emph{not subspaces} of $\sobh2(\W)$.  So before
      restricting, we need to extend and we use a traditional
      distribution-theory (or generalised-functions) approach.
      %%%%%%%%%%%%%%%%%%%%%%%%%%%%%%%%%%%%%%%%%%%%%%%%%%%%%%%%%%%%%%%%%%%%%%%%
      Given a function $v\in\sobh2(\W)$ and let $\geovec
      n:\partial\W\to\reals^d$ be the outward pointing normal of $\W$ then the
      Hessian of $v$, $\Hess v$ satisfies the following identity:
      \begin{equation}
        \label{eq:generalised-hessian}
        \ltwop{\Hess v}{\phi}
        = 
        -
        \int_\W \nabla v \otimes \nabla \phi
        +
        \int_{\partial\W} \nabla v \otimes \geovec n\,\phi 
        \Foreach \phi\in\sobh1(\W),
      \end{equation}
      where $\geovec a\otimes\geovec b:=\geovec a\transposevec b$ for
      $\vec a,\vec b$ column vectors in $\R d$.  If $v\in\sobh{1}(\W)$
      with $\restriction{\grad v}{\boundary\W}\in\sobh{-1/2}(\partial\W)$ the
      right-hand side of (\ref{eq:generalised-hessian}) still makes sense
      and defines $\Hess v$ as an element in the dual of $\sobh1(\W)$ via
      %%.  In this case we may define a
      %%functional, $\H[v]$, such that
      \begin{equation}
        %%\ltwop{\H[v]}{\phi} 
        %%= 
        \duality{\Hess v}{\phi}
        :=
        -
        \int_\W \nabla v \otimes \nabla \phi
        +
        \int_{\partial\W} \nabla v \otimes \geovec n\,\phi
        \Foreach \phi\in\sobh1(\W)
        ,
      \end{equation}
      where $\duality\cdot\cdot$ denotes the duality action on $\sobh1(\W)$
      from its dual.  We call $\Hess v$ the \emph{generalised Hessian} of
      $v$, and assuming that the coefficient tensor $\A$ is in
      $\cont0(\W)^{d\times d}$, for the product with a distribution to make
      sense, we now seek $u\in\sobhz1(\W)$ such that $\restriction{\grad
        u}{\boundary\W}\in\sobh{-1/2}(\W)$ and whose generalised Hessian
      satisfies
      \begin{equation}
        \label{eq:linear-model}
        \duality{\frob\A\Hess v}{\phi}
        =
        \ltwop f\phi
        \Foreach
        \phi\in\sobh1(\W).
      \end{equation}
      %%%%%%%%%%%%%%%%%%%%%%%%%%%%%%%%%%%%%%%%%%%%%%%%%%%%%%%%%%%%%%%%%%%%%%%%
      \begin{comment}
        %% I don't think this is correct: the operator \H works only at the
        %% discrete level it should be called \H_\fes really...
        %%%%%%%%%%%%%%%%%%%%%%%%%%%%%%%%%%%%%%%%%%%%%%%%%%%%%%%%%%%%%%%%%%%%%%%%
        The mixed formulation of the model problem (\ref{Problem}) we consider is to seek
        the pair $(u, \H[u]) \in \sobh{1}(\W)\cap\sobh{1}(\partial\W)
        \times \leb{2}(\W)^{d\times d}$ such that
        \begin{gather}
          \label{eq:mixed-1}
          \ltwop{\H[u]}{\phi} 
          +
          \int_\W{\nabla u}\otimes{\nabla \phi} 
          -
          \int_{\partial \W} {\nabla u}\otimes{\geovec{n} \ \phi}
          =
          \geomat 0
          \\
          \label{eq:mixed-2}
          \ltwop{\frob{\A}{\H[u]}}{\psi} 
          =
          \ltwop{f}{\psi}
          \Foreach (\phi,\psi)\in\sobh1(\W)\times \leb{2}(\W).
        \end{gather}
        %%%%%%%%%%%%%%%%%%%%%%%%%%%%%%%%%%%%%%%%%%%%%%%%%%%%%%%%%%%%%%%%%%%%%%%%
      \end{comment}
      %%%%%%%%%%%%%%%%%%%%%%%%%%%%%%%%%%%%%%%%%%%%%%%%%%%%%%%%%%%%%%%%%%%%%%%%
      \subsection{Finite element discretisation and finite element Hessian}
      We discretise (\ref{eq:linear-model}) for simplicity
      with a standard piecewise polynomial approximation for test and trial
      spaces for both problem variable, $U$, and auxiliary (mixed-type) variable,
      $\H[U]$. Let $\T{}$ be a conforming, shape regular
      triangulation of $\W$, namely, $\T{}$ is a finite family of sets such
      that
      \begin{enumerate}
      \item $K\in\T{}$ implies $K$ is an open simplex (segment for $d=1$,
        triangle for $d=2$, tetrahedron for $d=3$),
      \item for any $K,J\in\T{}$ we have that $\closure K\meet\closure J$ is
        a full subsimplex (i.e., it is either $\emptyset$, a vertex, an
        edge, a face, or the whole of $\closure K$ and $\closure J$) of both
        $\closure K$ and $\closure J$ and
      \item $\union{K\in\T{}}\closure K=\closure\W$.
      \end{enumerate}
      We use the convention where $\funk h\W\reals$ denotes the
      \emph{meshsize function} of $\T{}$, i.e.,
      \begin{equation}
        h(\vec{x}):=\max_{\closure K\ni \vec x}h_K.
      \end{equation}
      We introduce the \emph{finite element spaces}
      \begin{gather}
        \label{eqn:def:finite-element-space}
        \fes
        :=\ensemble{\Phi \in \sobh1(\W)}{\Phi\vert_{K} \in \poly p\Forall
          K\in\T{} \AND \Phi\in\cont{0}(\W)},
        \\
        \feszero
        :=\fes \cap \hoz,
      \end{gather}
      where $\poly k$ denotes the linear space of polynomials in $d$
      variables of degree no higher than a positive integer $k$. We
      consider $p\geq 1$ to be fixed and denote by $\Nzero :=
      \dim{\feszero}$ and $N := \dim{\fes}$. 
      %%%%%%%%%%%%%%%%%%%%%%%%%%%%%%%%%%%%%%%%%%%%%%%%%%%%%%%%%%%%%%%%%%%%%%%%
      %%\subsection{The nonvariational FE discretisation of problem (\ref{Problem})}
      The discretisation of problem then reads: Find $(U, \H[U]) \in
      \feszero \times \fes^{d\times d}$
      such that
      \begin{equation}
        \begin{gathered}
          \label{eq:discrete-linear}
          \ltwop{\H[U]}{\Phi} 
          = 
          -
          \int_\W{\nabla U}\otimes{\nabla \Phi} 
          +
          \int_{\partial \W} {\nabla U}\otimes{\geovec{n} \ \Phi}
          \Foreach \Phi\in\fes,
          \\
          \ltwop{\frob{\A}{\H[U]}}{\Psi} 
          =
          \ltwop{f}{\Psi}
          \Foreach \Psi\in\feszero.
        \end{gathered}
      \end{equation}
      For an algebraic formulation of (\ref{eq:discrete-linear}) we refer
      the reader to \cite[\S 2]{LakkisPryer:2011}.  Note that this
      discretisation can be interpreted as a mixed method whereby the
      first (matrix) equation defines the \emph{finite element Hessian} and the
      second (scalar) equation approximates the original PDE (\ref{Problem}).
      %%%%%%%%%%%%%%%%%%%%%%%%%%%%%%%%%%%%%%%%%%%%%%%%%%%%%%%%%%%%%%%%%%%%%%%%
      \subsection{Two discretisation stategies of \eqref{eq:modelproblem}}
        \label{sec:two-nonlinear-finite-element-methods}
        The finite element Hessian allows us two discretisation strategies.
        The first strategy, detailed in \S\ref{sec:unconstrained-fnl-pde},
        consists in applying Newton first to set-up
        (\ref{eq:system-of-nonvar-pdes}) and then using the NVFEM
        (\ref{eq:discrete-linear}) to solve each step.  A second strategy
        becomes possible, upon noting that given $U\in\fes$ the finite
        element Hessian $\H[U]$ is a regular function,\footnote{A
          generalised function $v$ is a \emph{regular function}, or just
          \emph{regular}, if it can be represented by a Lebesgue measurable
          function $f\in\lebloc1$ such that $\duality v\phi=\int_\W f\phi$
          for all $\phi\in\cont\infty_0(\W)$.  We follow the customary and
          harmless abuse in identifying $v$ with $f$.}  which the
        generalised Hessian $\Hess U$ might fail to be.  This allows to
        apply nonlinear functions such as $F$ to $\H[U]$ and consider the
        following \emph{fully nonlinear finite element method} (FNFEM)
        \begin{equation}
        \begin{gathered}
          \label{eq:fully-nonlinear-FEM}
          \ltwop{\H[U]}{\Phi} 
          = 
          -
          \int_\W{\nabla U}\otimes{\nabla \Phi} 
          +
          \int_{\partial \W} {\nabla U}\otimes{\geovec{n} \ \Phi}
          \Foreach \Phi\in\fes,
          \\
          \ltwop{F(\H[U])}\Psi
          =
          \ltwop f\Psi
          \Foreach \Psi\in\fez.
        \end{gathered}
        \end{equation}
        Of course, in order to solve the second equation, a
        finite-dimensional Newton method may be necessary (but this
        strategy leaves the door open for other nonlinear solvers, e.g.,
        fixed point iterations).  A finite element code based on this idea
        will be tested in \S\ref{sec:Pucci} to solve the Pucci equation.
      
        In summary the finite element Hessian allows both paths in the
        following diagram:
        \begin{equation}
          \small
          \label{eqn:dia:Newton-discretizations}
          \begin{tikzpicture}
            \path (0,3) node(X){fully nonlinear PDE \eqref{eq:modelproblem}};
            \path (6,3) node(N){nonvariational linear PDE's \eqref{eq:system-of-nonvar-pdes}};
            \path (0,0) node(F){fully nonlinear FE discretization \eqref{eq:fully-nonlinear-FEM}};
            \path (6,.5) node(D){\color adiscrete linear 1};
            \path (6,0) node(L){\color bdiscrete linear 2};
            \draw[-stealth,color=a] (X)--(N) node[pos=0.5,above]{Newton};
            \draw[-stealth,color=a] (N)--(D) node[pos=0.5,above,sloped]{NVFEM};
            \draw[-stealth,color=b] (X)--(F) node[pos=0.5,above,sloped]{FNFEM};
            \draw[-stealth,color=b] (F)--(L) node[pos=0.5,above]{Newton};
          \end{tikzpicture}
        \end{equation}
        Although the diagram in (\ref{eqn:dia:Newton-discretizations}) does
        not generally commute, if the function $F$ is algebraically
        accessible, then it is commutative.  By ``algebraically accessible''
        we mean a function that can be computed in a finite number of
        algebraic operations or inverses thereof.  In this paper, we use
        only algebraically accessible nonlinearities, but, in principle
        assuming derivatives are available, our methods could be extended to
        algebraically inaccessible nonlinearities, such as Bellman's (or
        Isaacs's) operators involving optimums over infinite families, e.g.,
        \begin{equation}
          F(\geomat M):=\inf_{\alpha\in\cA}\frob{\geomat L_{\alpha}}{\geomat M}
          \quad
          \qpbig{\OR
          F(\geomat M):=\inf_{\alpha\in\cA}\sup_{\beta\in\cB}\frob{\geomat L_{\alpha,\beta}}{\geomat M}}
          ,
        \end{equation}
        where $\ensemble{\geomat L_{\alpha}}{\alpha\in\cA}$ 
        (or $\ensemble{\geomat L_{\alpha,\beta}}{(\alpha,\beta)\in\cA\times\cB}$)
        is a family of elliptic operators.
      %%%%%%%%%%%%%%%%%%%%%%%%%%%%%%%%%%%%%%%%%%%%%%%%%%%%%%%%%%%%%%%%%%%%%%%%
      %%%%%%%%%%%%%%%%%%%%%%%%%%%%%%%%%%%%%%%%%%%%%%%%%%%%%%%%%%%%%%%%%%%%%%%%
      \section{The discretisation of unconstrained fully nonlinear PDEs}
      \label{sec:unconstrained-fnl-pde}
      
      In this section we detail the application of the method reviewed in
      \S\ref{sec:ndfem} to the fully nonlinear model problem
      \eqref{eq:modelproblem}.  Many fully nonlinear elliptic PDEs must be
      constrained in order to admit a unique solution. For example the \MAD
      is elliptic and admits a unique solution in the cone of convex (or
      concave) functions when $f>0$ (or $f<0$, respectively).  Before we
      turn our attention to the more complicated constrained PDE's in
      \S\ref{sec:ma} and we illustrate the Newton--NVFEM method in the simplest light.
      In this section we study fully nonlinear PDEs which have no such
      constraint.
      %%%%%%%%%%%%%%%%%%%%%%%%%%%%%%%%%%%%%%%%%%%%%%%%%%%%%%%%%%%%%%%%%%%%%%%%
      \begin{Hyp}[unconditionally elliptic linearisation]
        \label{ass:well-posed}
        We assume, in this section, that the Newton-step linearisation
        (\ref{eq:system-of-nonvar-pdes}) is elliptic.  For this assumption
        to hold, it is sufficient to assume uniform ellipticity, i.e.,
        (\ref{eq:ellipticity-for-derivative}) with $\cC=\symm$ and $\mu>0$
        independent of $\geomat M$.
      \end{Hyp}
      %%%%%%%%%%%%%%%%%%%%%%%%%%%%%%%%%%%%%%%%%%%%%%%%%%%%%%%%%%%%%%%%%%%%%%%%
      \subsection{The Newton-NVFEM method}
        \label{The:nlfem}
        Suppose we are given a BVP of the form, finding
        $u\in\sobh2(\W)\cap\hoz$ such that
        \begin{equation}
          \label{eq:nonlinearproblem}
          \begin{split}
            \nlop{u} = \nlfunk(\Hess u) - f &= 0 \qquad \text{ in } \W,     
          \end{split}
        \end{equation}
        which satisfies Assumption \ref{ass:well-posed}.
      
        Upon applying Newton's method to approximate the solution of problem
        (\ref{eq:nonlinearproblem}) we obtain a sequence of functions
        $(u^n)_{n\in\naturals_0}$ solving the following linear equations in
        nonvariational form,
        \begin{equation}
          \label{eq:linearisedequation}
          \frob{\N(\Hess u^n)}{\Hess u^{n+1}} = g(\Hess u^n)
        \end{equation}
        where
        \renewcommand{\thefootnote}{\fnsymbol{footnote}} 	
        \begin{gather}
          \label{eq:newton_prob_data:N}
          \N(\geomat X) := F'(\geomat X),\\
          \label{eq:newton_prob_data:g}
          g(\geomat X) 
          :=
          f - \nlfunk(\geomat X) + \frob{F'(\geomat X)}{\geomat X}.
        \end{gather}
      
        The nonlinear finite element method to approximate
        (\ref{eq:linearisedequation}) is: given an initial guess $U^0 :=
        \Pi_0 u^0$ for each $n\in\naturals_0$ find $(U^{n+1},
        \H[U^{n+1}])\in \feszero\times\fes^{d\times d}$ such that
        \begin{equation}
          \begin{gathered}
            \label{eq:discretenewtonsmethod}
            \ltwop{\H[U^{n+1}]}{\Phi} 
            +
            \int_\W{\nabla U^{n+1}}\otimes{\nabla \Phi} 
            -
            \int_{\partial \W} {\nabla U^{n+1}}\otimes{\geovec{n}}\,\Phi
            =
            \geomat 0
            \Foreach \Phi\in\fes
            \\
            \AND
            \ltwop{\frob{\N(\H[U^{n}])}{\H[U^{n+1}]}}{\Psi} 
            =
            \ltwop{g(\H[U^n])}{\Psi}
            \Foreach \Psi\in\feszero.
          \end{gathered}
        \end{equation}
      %%%%%%%%%%%%%%%%%%%%%%%%%%%%%%%%%%%%%%%%%%%%%%%%%%%%%%%%%%%%%%%%%%%%%%%%
        \begin{Obs}[initial conditions]
          
        \end{Obs}
      \subsection{Numerical experiments: a simple example}
\begin{figure}[h]
        \caption[]{\label{Fig:nonlin-regular-abs-lap} Numerical experiments
          for Example \ref{ex:fullynl-abs}. Choosing $f$ appropriately such
          that $u(\geovec x) = \exp\qp{-10\norm{\geovec x}^2}$. We use an
          initial guess $u^0 = 0$ and run the iterative procedure until
          $\Norm{U^{n+1} - U^n} \leq 10^{-8}$, setting $U := U^M$
          the final Newton iterate of the sequence. Here we are plotting
          log--log error plots together with experimental convergence rates
          for $\leb{2}(\W), \sobh1(\W)$ error functionals for the problem
          variable, $U$, and an $\leb{2}(\W)$ error functional for the
          auxiliary variable, $\H[U]$. Notice that there is a
          ``superconvergence'' of the auxiliary variable for both
          approximations.}
        % \begin{center} 
        \subfigure[][{Taking $\fes$ to be the space of
            piecewise linear functions on $\W$ ($p = 1$). Notice that $\Norm{u -
              U^M} = \Oh(h^2)$, $\norm{u - U^M}_1 = \Oh(h)$ and $\Norm{\Hess u
              - \H[U^M]} = \Oh(h^{0.5})$. }]{
          \label{Fig:sin-lap-p1}
          \includegraphics[scale=\figscale,width=0.47\figwidth]{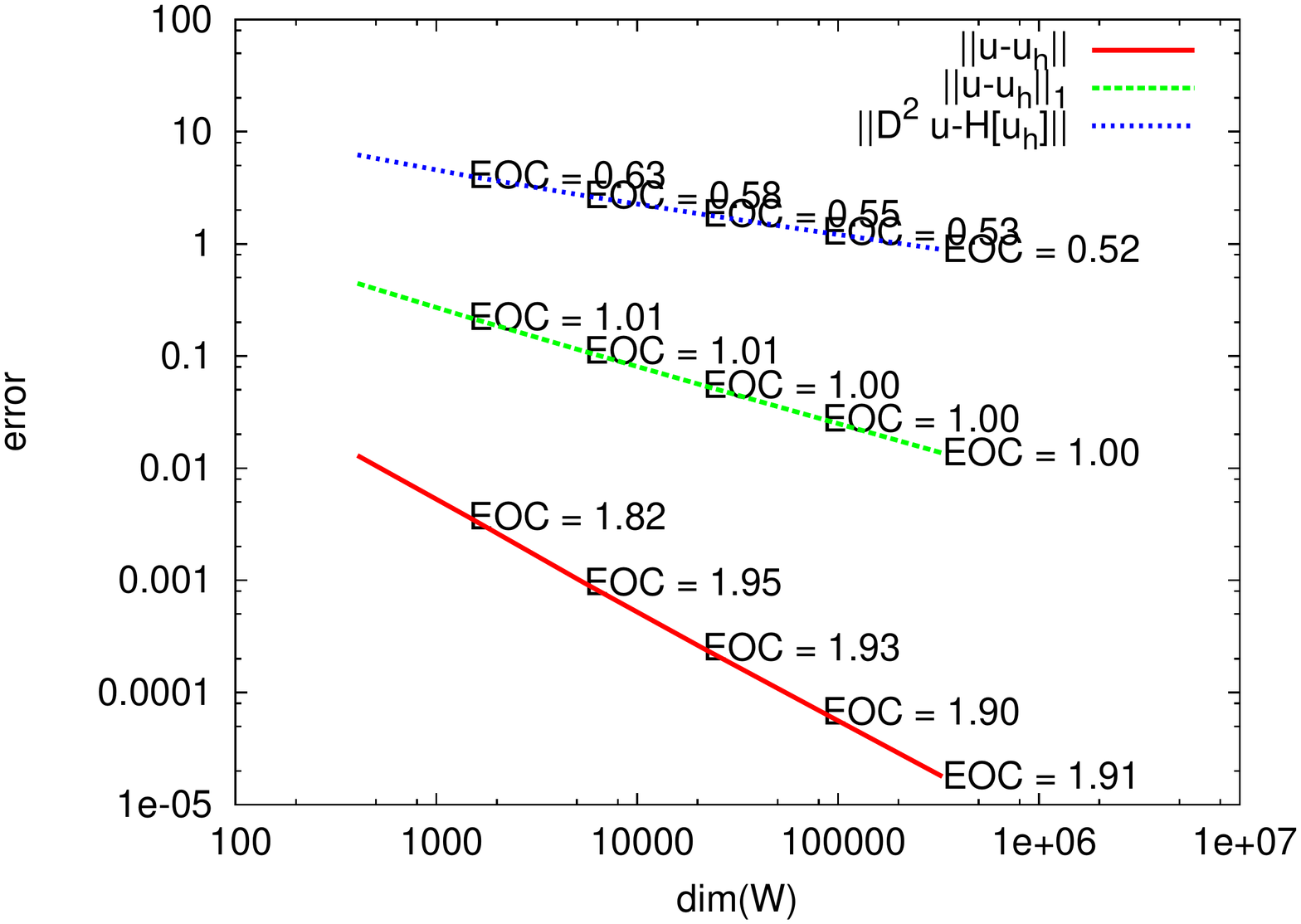}
        } \hfill \subfigure[][{Taking $\fes$ to be the space of piecewise
            quadratic functions on $\W$ ($p = 2$). Notice that $\Norm{u -
              U^M} = \Oh(h^3)$, $\norm{u - U^M}_1 = \Oh(h^2)$ and
            $\Norm{\Hess u - \H[U^M]} = \Oh(h^{1.5})$}]{
          \label{Fig:sin-lap-p2}
          \includegraphics[scale=\figscale,width=0.47\figwidth]{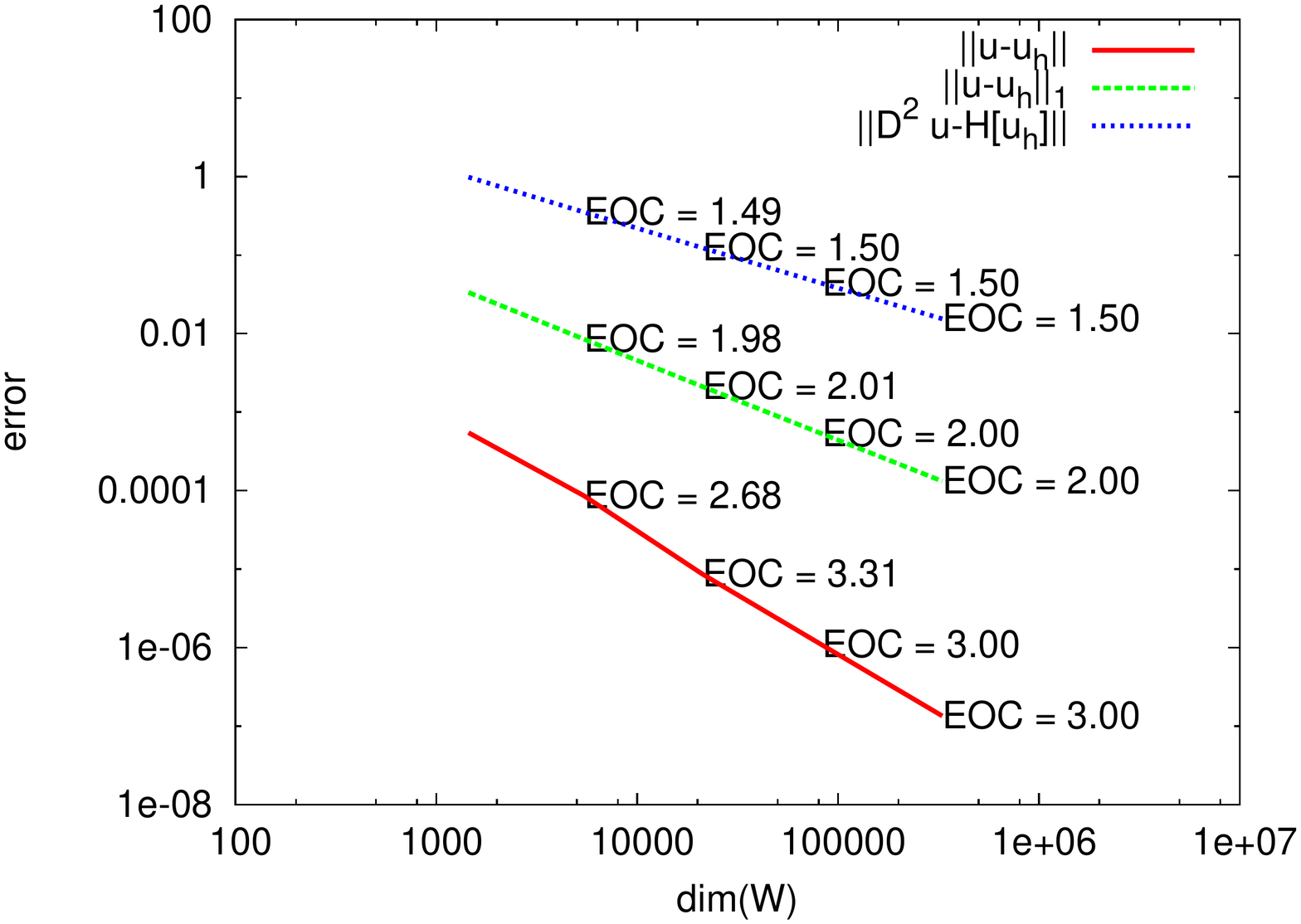}
        }
        %\end{center}
      \end{figure}

      %%%%%%%%%%%%%%%%%%%%%%%%%%%%%%%%%%%%%%%%%%%%%%%%%%%%%%%%%%%%%%%%%%%%%%%%
      In this section we detail numerical experiments aimed at demonstrating
      the application of \eqref{eq:discretenewtonsmethod} to a simple model
      problem.
      %%%%%%%%%%%%%%%%%%%%%%%%%%%%%%%%%%%%%%%%%%%%%%%%%%%%%%%%%%%%%%%%%%%%%%%%
      \begin{Example}[a simple fully nonlinear PDE]
        \label{ex:fullynl-abs}
        The first example we consider is a fully nonlinear PDE with a very
        smooth nonlinearity. The problem is
        \begin{equation}
            \begin{split}
            \nlop{u} := \sin{\Delta u} + 2\Delta u -f &= 0 \text{ in } \W,
            \\
            u &= 0 \text{ on } \pd{}\W.
            \end{split}
            %%\begin{split}
            %%       \nlop{u} := \norm{\Delta u} + 2\Delta u - f &= 0 \quad \text{ in } \W
            %%       \\
            %%       u &= 0 \quad \text{ on } \partial\W
            %%\end{split}
        \end{equation}
        which is specifically constructed to be uniformly elliptic. Indeed
        \begin{equation}
          \nlfunk'(\Hess u)
          = 
          \left(\cos{\Delta u} + 2\right)\eye.
        \end{equation}
        which is uniformly positive definite.
        The Newton linearisation of the problem is then: Given $u^0$, for
        $n\in\naturals_0$ find $u^{n+1}$ such that
        \begin{equation}
          \frob{\left(\cos{\Delta u^n} + 2\right)\eye}
               {\Hess (u^{n+1} - u^n)}
               =
               f - \sin{\Delta u^n} - 2\Delta u^n.
        \end{equation}
        and our approximation scheme is nothing but
        \ref{eq:discretenewtonsmethod} with
        \begin{gather}
          \N(\geomat X) = \left(\cos{\trace{\geomat X}} + 2\right)\geomat I
          \\
          g(\geomat X) = f - \sin{\trace{\geomat X}} - 2\trace{\geomat X}.
        \end{gather}
        Figure \ref{Fig:nonlin-regular-abs-lap} details a numerical
        experiment on this problem when $d=2$ and when $\W =[-1,1]^2$ is a
        square which is triangulated using a criss-cross mesh.
      \end{Example}
      %%%%%%%%%%%%%%%%%%%%%%%%%%%%%%%%%%%%%%%%%%%%%%%%%%%%%%%%%%%%%%%%%%%%%%%%
      \begin{Rem}[simplification of Example \ref{ex:fullynl-abs}]
        \label{rem:variational-simplification}
        Example \ref{ex:fullynl-abs} can be simplified considerably by noticing that
        \begin{equation}
          \int_\W \tr{\H[U]} \Phi = \int_\W \Transpose{\qp{\nabla U}}\nabla \Phi
          \Foreach \Phi\in\feszero.
        \end{equation}
        This coincides with the definition of the \emph{discrete Laplacian}
        and makes the NVFEM coincide with the standard conforming FEM.  This
        observation applies to all fully nonlinear equations, with
        nonlinearity of the form \eqref{eq:modelproblem} with
        \begin{equation}
          F(\geomat M):=a(\tr{\geomat M}),
        \end{equation}
        for some given $a$.  This class of problems, can be solved using a
        variational finite element method and can be used for comparison
        with our method. Note that in
        \cite{JensenSmears:12:techreport:Finite} the authors use this
        together with a localisation argument in order to prove convergence
        of a finite element method for a specific class of
        Hamilton--Jacobi--Bellman equation.  Their method coincides with ours,
        for an appropriate choice of quadrature.
      \end{Rem}
      
      %%%%%%%%%%%%%%%%%%%%%%%%%%%%%%%%%%%%%%%%%%%%%%%%%%%%%%%%%%%%%%%%%%%%%%%%
      %%%%%%%%%%%%%%%%%%%%%%%%%%%%%%%%%%%%%%%%%%%%%%%%%%%%%%%%%%%%%%%%%%%%%%%%
      \begin{Example}[nonvariational example]
        \label{eq:fullynl-unsolvable}
        This is a simple example where the \emph{variational trick}
        mentioned in Remark \ref{rem:variational-simplification} cannot be applied. We fix
        $d=2$ and consider the problem
        \begin{equation}
          \begin{split}
            \nlop{u}:=
            \qp{\partial_{11} u}^3
            +
            \qp{\partial_{22} u}^3
            +
            \partial_{11} u
            +
            \partial_{22} u
            -
            f
            &=
            0 \text{ in } \W
            \\
            u &= 0 \text{ on } \partial\W.
          \end{split}
        \end{equation}
        The approximation scheme is then \eqref{eq:discretenewtonsmethod} with
        \begin{gather}
          \N(\geomat X) :=
          \begin{bmatrix}
            3\geomat X_{11}^2 + 1 & 0\\
            0 & 3\geomat X_{22}^2 +1
          \end{bmatrix}
          \\
          g(\geomat X)
          := 
          f + 2\qp{\geomat X_{11}^3 + \geomat X_{22}^3}.
        \end{gather}
        Figure \ref{Fig:nonlin-krylov} details a numerical experiment on
        this problem in the case $d=2$ and $\W = [-1,1]^2$ triangulated with
        a criss-cross mesh. A similar example is also studied in
        \cite[Ex 5.2]{DavydovSaeed:12:techreport} using B\"ohmers method.
      
      %%%%%%%%%%%%%%%%%%%%%%%%%%%%%%%%%%%%%%%%%%%%%%%%%%%%%%%%%%%%%%%%%%%%%%%%
      \begin{figure}[h]
        \caption[]{\label{Fig:nonlin-krylov} Numerical experiments
          for Example \ref{eq:fullynl-unsolvable}. Choosing $f$ appropriately such
          that $u(\geovec x) = \exp\qp{-10\norm{\geovec x}^2}$. We use an
          initial guess $u^0 = 0$ and run the iterative procedure until
          $\Norm{U^{n+1} - U^n} \leq 10^{-8}$, setting $U := U^M$
          the final Newton iterate of the sequence. Here we are plotting
          log--log error plots together with experimental convergence rates
          for $\leb{2}(\W), \sobh1(\W)$ error functionals for the problem
          variable, $U$, and an $\leb{2}(\W)$ error functional for the
          auxiliary variable, $\H[U]$. Notice that there is a
          ``superconvergence'' of the auxiliary variable for both
          approximations.}
        % \begin{center} 
        \subfigure[][{Taking $\fes$ to be the space of
            piecewise linear functions on $\W$ ($p = 1$). Notice that $\Norm{u -
              U^M} = \Oh(h^2)$, $\norm{u - U^M}_1 = \Oh(h)$ and $\Norm{\Hess u
              - \H[U^M]} \approx \Oh(h^{1.5})$. }]{
          \label{Fig:sin-lap-p1}
          \includegraphics[scale=\figscale,width=0.47\figwidth]{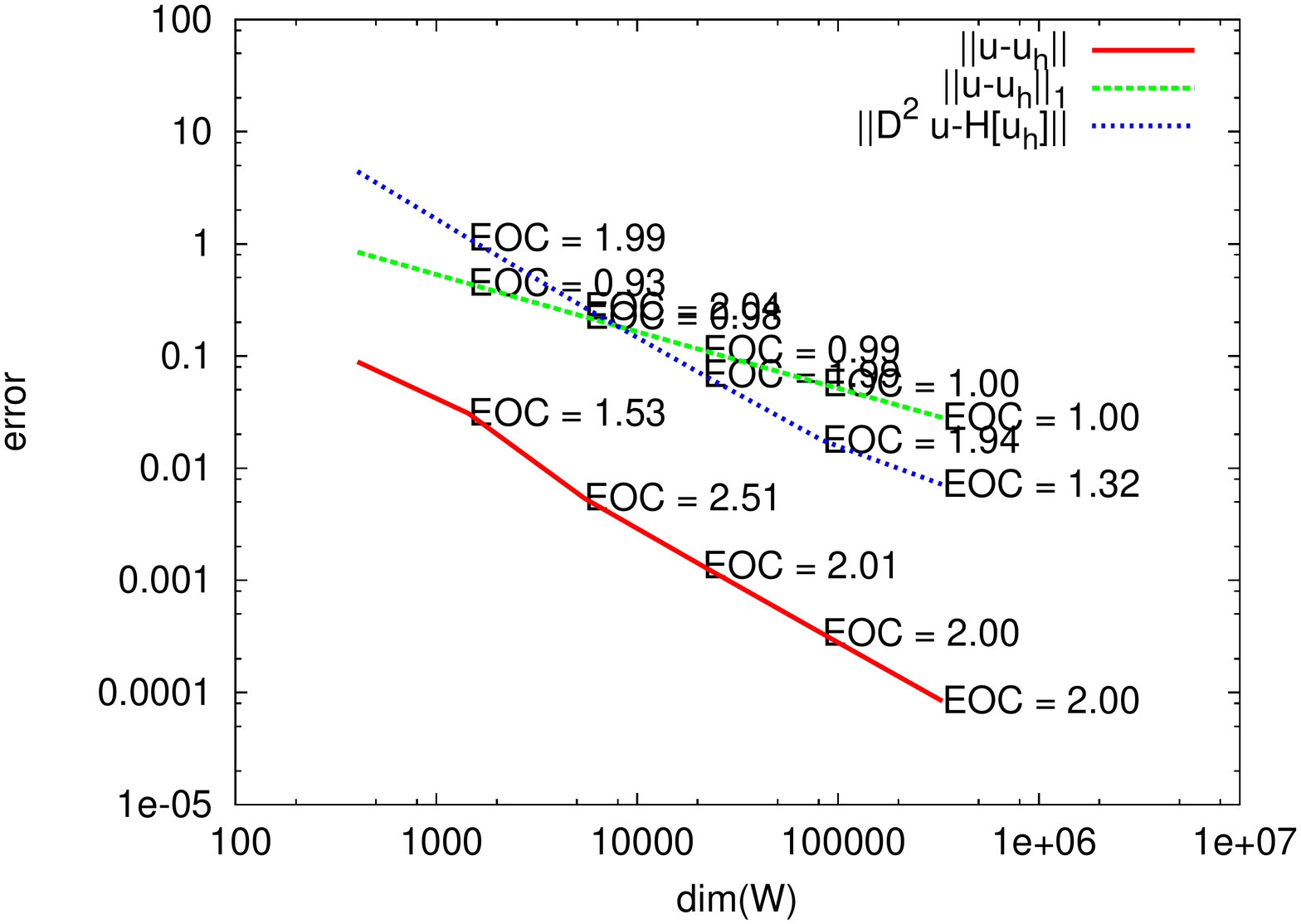}
        } \hfill \subfigure[][{Taking $\fes$ to be the space of piecewise
            quadratic functions on $\W$ ($p = 2$). Notice that $\Norm{u -
              U^M} = \Oh(h^3)$, $\norm{u - U^M}_1 = \Oh(h^2)$ and
            $\Norm{\Hess u - \H[U^M]} = \Oh(h^{1.3})$}]{
          \label{Fig:sin-lap-p2}
          \includegraphics[scale=\figscale,width=0.47\figwidth]{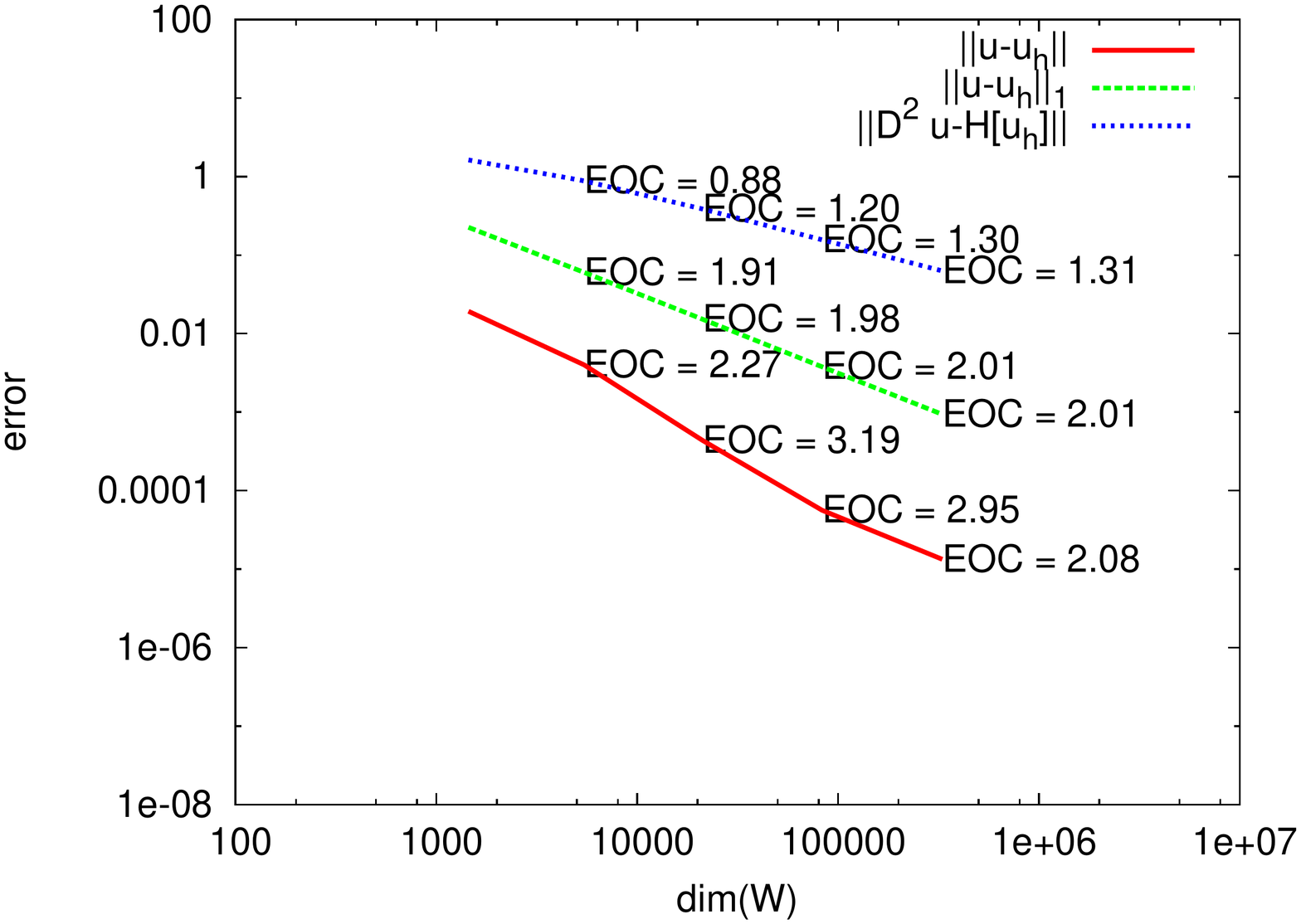}
        }
        %\end{center}
      \end{figure}
      %%%%%%%%%%%%%%%%%%%%%%%%%%%%%%%%%%%%%%%%%%%%%%%%%%%%%%%%%%%%%%%%%%%%%%%%
      
      \end{Example}
      
      %\margnote[omar]{Tristan, we should put back that example you had previosly with the cubes and all that.  Davydov and Saeed have included that in their paper on splines and we should keep it for comparison.}
      
      % Done now
%%%%%%%%%%%%%%%%%%%%%%%%%%%%%%%%%%%%%%%%%%%%%%%%%%%%%%%%%%%%%%%%%%%%%%%%
      %%%%%%%%%%%%%%%%%%%%%%%%%%%%%%%%%%%%%%%%%%%%%%%%%%%%%%%%%%%%%%%%%%%%%%%%
      %%% mode:latex ***
      %%% tex-main-file: "../nonlinear.tex"  ***
      %%%%%%%%%%%%%%%%%%%%%%%%%%%%%%%%%%%%%%%%%%%%%%%%%%%%%%%%%%%%%%%%%%%%%%%%
      \section{The \MAD problem}
      \label{sec:ma}
      
      In this section we propose a numerical method for the \MAD (MAD) problem
      \begin{equation}
        \label{eq:ma}
        \begin{split}
          \det{\Hess u}
          &=
          f \text{ in }\W
          \\
          u 
          &= 
          g \text{ on }\partial \W.
        \end{split}
      \end{equation}
      Our numerical experiments exhibit robustness of our method when
      computing (smooth) classical solutions of the MAD equation.  Most
      importantly we noted the following facts:
      \begin{enumerate}[(i)\ ]
      \item 
        the use of $\poly p$ elements with $p\geq2$ is essential
        as $\poly1$ do not work,
      \item 
        the convexity of the Newton iterates is
        conserved throughout the computation, in a similar way to the
        observations in \cite{loeper-rapetti:05}, where the authors prove this
        convexity-conservation property.
      \end{enumerate}
      Our observations are purely
      empirical from computations, which leaves an interesting open problem
      of proving this property.
      %%%%%%%%%%%%%%%%%%%%%%%%%%%%%%%%%%%%%%%%%%%%%%%%%%%%%%%%%%%%%%%%%%%%%%%%
      % is an important example of a fully nonlinear
      %elliptic PDE since it is used as a model for other fully nonlinear
      %PDEs. %It is derived from differential geometry. The problem is
      %\begin{equation}
      %  \label
      %  \begin{split}
      %    \nlop{u} := \det{\Hess u} - f &= 0 \text{ in } \W
      %    \\
      %    u &= 0 \text{ on } \pd{}\W.
      %  \end{split}
      %\end{equation}
      %This equation is clearly fully nonlinear for $d\geq 2$, in fact it is
      %multi-linear with respect to columns (or rows) of the Hessian. %This
      %makes for simpler computations when it comes to linearising the
      %problem.
      
      %There have been numerical studies on this problem in the context of
      %finite differences by Oliker and Prussner
      %\cite{OlikerPrussner:1988}. Feng and Neilan
      %\cite{FengNeilan:2007,FengNeilan:2008,FengMichael:2008} propose a
      %mixed finite element method using sequences of fourth order
      %quasi-linear equations. Despite these works this is still a tricky
      %problem to formulate correctly. There are restrictions that must be
      %put in place in order to guarantee ellipticity for example.
      
      \begin{Rem}[the MAD problem fails to satisfy Assumption \ref{ass:well-posed}]
        To clarify Assumption \ref{ass:well-posed} for the MAD problem
        (\ref{eq:ma}), in view of the characteristic expansion of
        determinant if $\geomat X,\geomat Y \in \symm$
        \begin{equation}
          \det\qp{\geomat X + \epsilon \geomat Y} 
          =
          \det{\geomat X} + \epsilon \frob{\cof\qp{\geomat X}}{\geomat Y} + \Oh(\epsilon^2),
        \end{equation}
        where $\cof{\geovec X}$ is the
        matrix of cofactors of $\geovec X$. Hence
        \begin{equation}
          \label{eq:der-of-det-is-cof}
          \nlfunk'(\geomat X) = \cof{\geomat X}
          .
        \end{equation}
        This implies that the linearisation of MAD is only well posed if we
        restrict the class of functions we consider to those $u$ that
        satisfy
        \begin{equation}
          \label{eqn:MAD-ellipticity-via-derivative}
          \Transpose{\geovec \xi} \cof \Hess u\, \geovec \xi 
          \geq
          \lambda \norm{\geovec\xi}^2 \Foreach \geovec\xi\in\reals^d
        \end{equation}
        for some ($u$-dependent) $\lambda>0$.
        Note that (\ref{eqn:MAD-ellipticity-via-derivative}) is equivalent to the
        following two conditions as well
        \begin{equation}
          \begin{gathered}
            \Transpose{\geovec \xi} \Hess u\, \geovec \xi 
            \geq
            \lambda \norm{\geovec\xi}^2 \Foreach \geovec\xi\in\reals^d
            \\
            $u$ \text{ is strictly convex. }
          \end{gathered}
        \end{equation}
      
        \cite{LoeperRapetti:2005} have shown that for the \emph{continuous}
        (infinite dimensional) Newton method described in
        \ref{sec:newtons-method}, given an strictly convex initial guess
        $u^0$, each iterate $u^n$ will be convex.  It is crucial that this
        property is preserved at the discrete level, as it guarantees the
        solvability of each iteration in the \emph{discretised} Newton
        method.  For this it the right notion of convexity turns out to be
        the \emph{finite element convexity} as developed in
        \cite{Aguilera:2008}.  In \cite{Pryer:2010}, an intricate method
        based on semidefinite programming provided a way to constrain the
        solution in the case of $\poly1$ elements.  Here we observe that
        the finite element convexity is automatically preserved, provided we
        use $\poly2$ or higher conforming elements.
      \end{Rem}
      %%%%%%%%%%%%%%%%%%%%%%%%%%%%%%%%%%%%%%%%%%%%%%%%%%%%%%%%%%%%%%%%%%%%%%%%
      %We require $\W$ to be convex and $f > 0$ for a classical solution to
      %even exist. \MA (\ref{eq:ma}) will be uniformly elliptic if $\Hess u$
      %is positive definite. If these restrictions are then satisfied then
      %(\ref{eq:ma}) admits a unique convex viscosity solution.
      %%%%%%%%%%%%%%%%%%%%%%%%%%%%%%%%%%%%%%%%%%%%%%%%%%%%%%%%%%%%%%%%%%%%%%%%
      %Indeed without the constraint on the Hessian of the solution the
      %problem admits two solutions, one convex and one concave, in $d=2$ for
      %example with homogeneous Dirichlet boundary both $u$ and $-u$ solve
      %(\ref{eq:ma}).
      %%%%%%%%%%%%%%%%%%%%%%%%%%%%%%%%%%%%%%%%%%%%%%%%%%%%%%%%%%%%%%%%%%%%%%%%
      \subsection{Newton's method applied to \MA}
      
      In view of (\ref{eq:der-of-det-is-cof}) it is clear that
      \begin{equation}
        \Dnlop{u}v = \frob{\cof{\Hess u}}{\Hess v}.
      \end{equation} 
      Applying the methodology set out in \S\ref{sec:unconstrained-fnl-pde}
      we set
      \begin{gather}
        \label{eq:maN}
        \N(\Hess u^n) = \cof{\Hess u^n},
        \\
        \label{eq:mag}
        g(\Hess u^n) = f - \det{\Hess u^n} + \frob{\cof{\Hess
            u^n}}{\Hess u^n},
      \end{gather}
      
      \begin{Rem}[{relating cofactors to determinants}]
        \label{rem:cof-to-det}
        For a generic (twice differentiable) function $v$ it holds
        that
        \begin{equation}
          d \det{\Hess v} = \frob{\cof{\Hess v}}{\Hess v}.
        \end{equation}
        Using this formulation we could construct a simple fixed point
        method for the \MA equation.
      \end{Rem}
      In view of Remark \ref{rem:cof-to-det} $g$ can be further
      simplified
      \begin{equation}
        \begin{split}
          g(\Hess u^n) &= f - \det{\Hess u^n} + \frob{\cof{\Hess
              u^n}}{\Hess u^n}
          \\
          &= f + (d-1) \det{\Hess u^n}.
        \end{split}
      \end{equation}
      Newton's method reads: Given $u^0$ for each $n\in\naturals_0$ find
      $u^{n+1}$ such that
      \begin{equation}
        \label{eq:linearised-ma}
        \frob{\N(\Hess u^n)}{\Hess u^{n+1}} 
        =
        g(\Hess u^n).
      \end{equation}
      
      \subsection{Numerical experiments}
      
      In this section we study the numerical behaviour of the scheme
      presented in Definition \ref{The:nlfem} applied to the MAD problem.
      
      We present a set of benchmark problems constructed from the problem
      data such that the solution to the \MA equation is known. We fix $\W$
      to be the square $S = [-1,1]^2$ or $[0,1]^2$ (specified in the
      problem) and test convergence rates of the discrete solution to the
      exact solution.
      
\begin{figure}[h]
        \caption{ \label{Fig:EOC-ma1} Numerical results for the MAD problem
          on the square $S = [-1,1]^2$. We choose the problem data $f$ and
          $g$ appropriately such that the solution is the radially symmetric
          function $u(\geovec x) = \exp\qp{{\norm{\geovec x}^2}/{2}}$. We
          plot the finite element solution together with a log--log error
          plot for various error functionals as in Figure
          \ref{Fig:nonlin-regular-abs-lap}. Note for $p=2$ the $\leb{2}(\W)$
          error rate of convergence is suboptimal, this is in agreement with
          the numerical examples produced in \cite{BrennerGudiNeilanSung:2011}}
        \begin{center}
              \subfigure[][{The FE approximation to the function $u(\geovec x) = \exp\qp{\frac{\norm{\geovec
                x}^2}{2}}$.}]{
      %      \label{Fig:ma-radial}
        \includegraphics[scale=\figscale,width=0.47\figwidth]{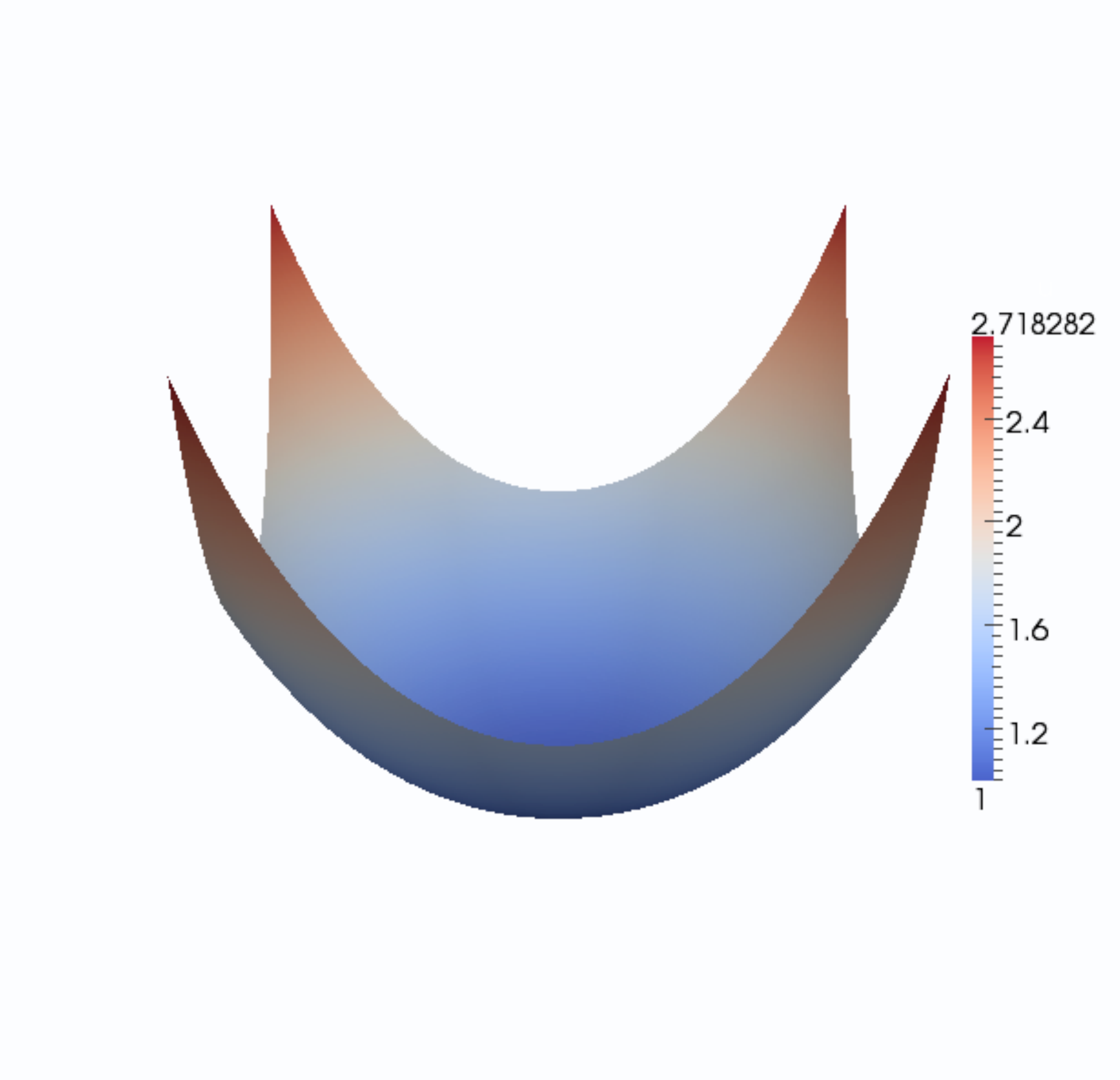}
      }  \hfill
          \subfigure[][{Log--log error plot for $\poly{2}$ Lagrange FEs.}]{
            \label{Fig:ma-radial}
            \includegraphics[scale=\figscale,width=0.47\figwidth]{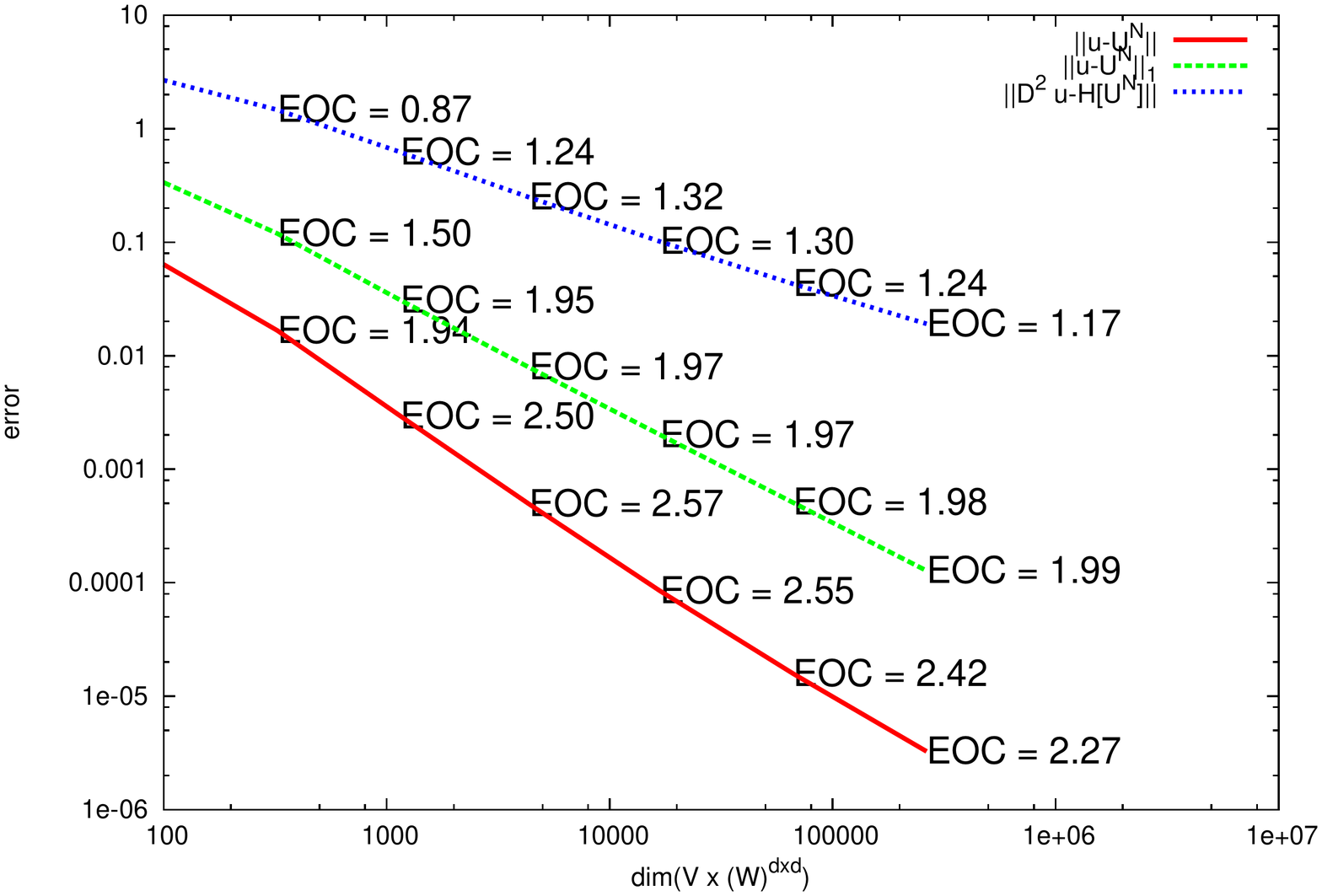}
          } 
        \end{center}
      \end{figure}

\begin{figure}[h]
        \caption{ \label{Fig:EOC-ma3} A FE--convexity test for the numerical example given in \ref{Fig:EOC-ma1}. We plot $\det{\H[U]}$ together with the principle minor of $\H[U]$.}
        \begin{center}
              \subfigure[][{The principal minor of $\H[U]$, an approximation to the Hessian of the function $u(\geovec x) = \exp\qp{\frac{\norm{\geovec x}^2}{2}}$.}]{
                \includegraphics[scale=\figscale,width=0.47\figwidth]{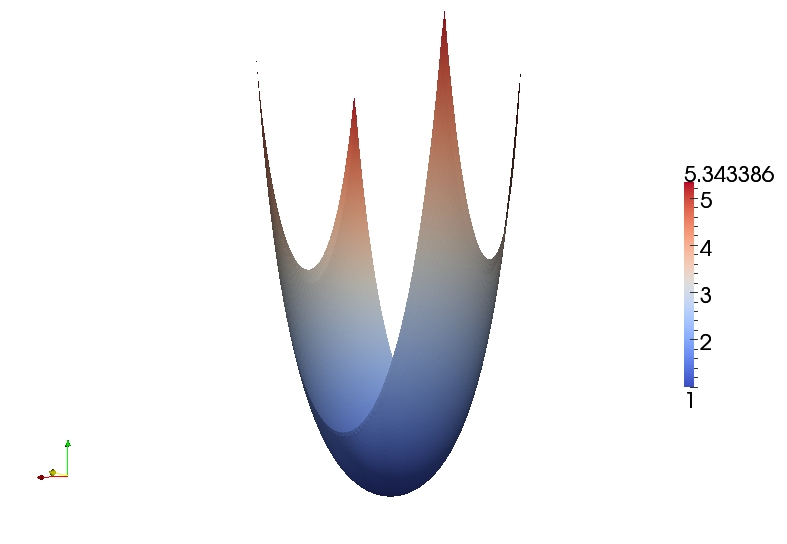}
      }  \hfill
          \subfigure[][{The determinant of $\H[U]$.}]{
            \includegraphics[scale=\figscale,width=0.47\figwidth]{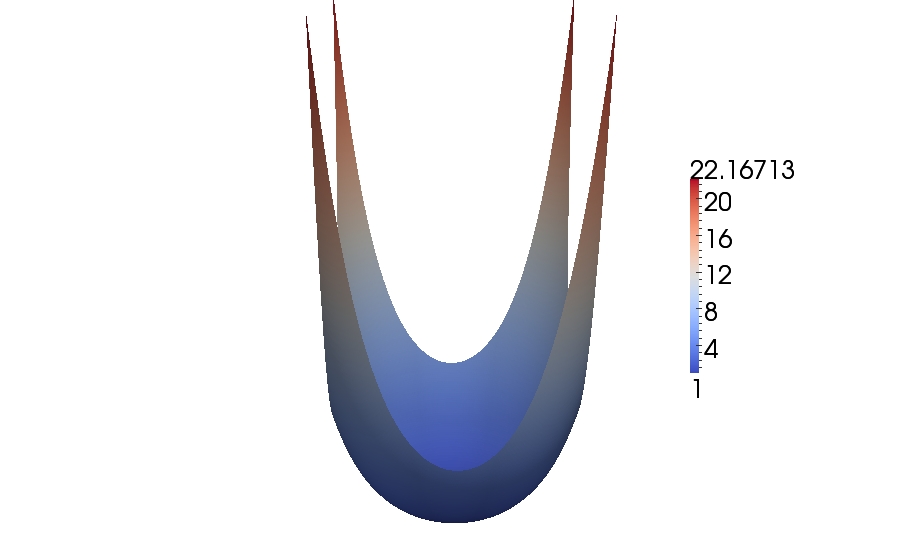}
          } 
        \end{center}
      \end{figure}

\begin{figure}[h]
        \caption{ \label{Fig:EOC-ma4} Numerical results for the MAD problem
          on the square $S = [-1,1]^2$. We choose the problem data $f$ and
          $g$ appropriately such that the solution is the radially symmetric
          function $u(\geovec x) = \exp\qp{{\norm{\geovec x}^2}/{2}}$.}
        \begin{center}
              \subfigure[][{The $\poly{1}$ FE approximation to the function $u(\geovec x) = \exp\qp{\frac{\norm{\geovec
                x}^2}{2}}$.}]{
                \includegraphics[scale=\figscale,width=0.47\figwidth]{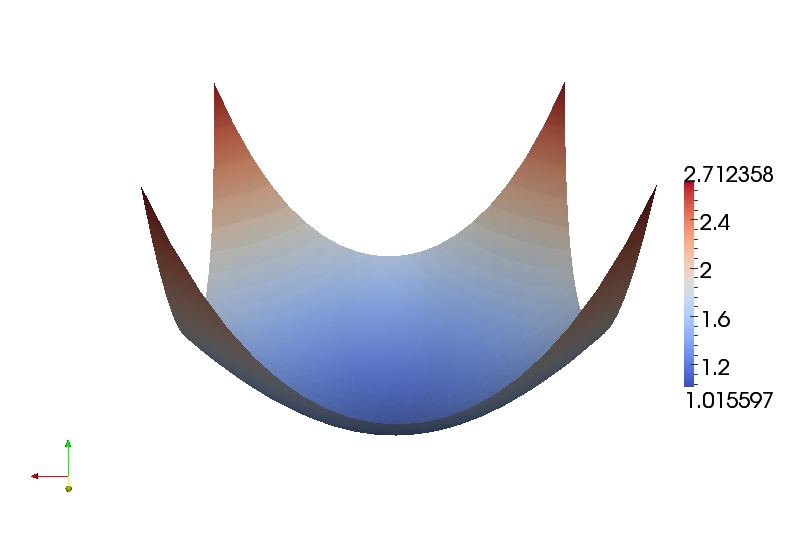}
      }  \hfill
          \subfigure[][{The error $u-U$ plotted as a function over $\W$. Note the 
      FE approximation does not converge in this case, see Remark \ref{rem:degree-of-fe-space}}]{
            \includegraphics[scale=\figscale,width=0.47\figwidth]{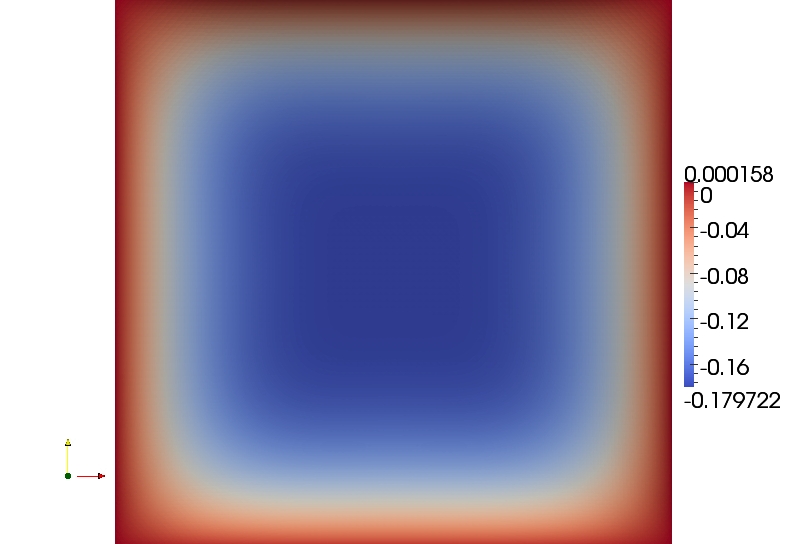}
          } 
        \end{center}
      \end{figure}

\begin{figure}[h]
        \caption{\label{Fig:EOC-ma2} Numerical results for the MAD problem
          on the square $S=[-1,1]^2$. Choosing $f$ and $g$ appropriately
          such that the solution is $u(\geovec x) = -\qp{2 - x_1^2 -
            x_2^2}^{1/2}$. Note the function has singular derivatives on the
          corners of $S$. We plot the finite element solution together with
          a log--log error plot for various error functionals as in Figure
          \ref{Fig:nonlin-regular-abs-lap}.}
        \begin{center}
          \subfigure[][{The FE approximation to the function $u(\geovec x)
              = -\qp{2 - x_1^2 - x_2^2}^{1/2}$.}]{
      %      \label{Fig:ma-radial}
        \includegraphics[scale=\figscale,width=0.47\figwidth]{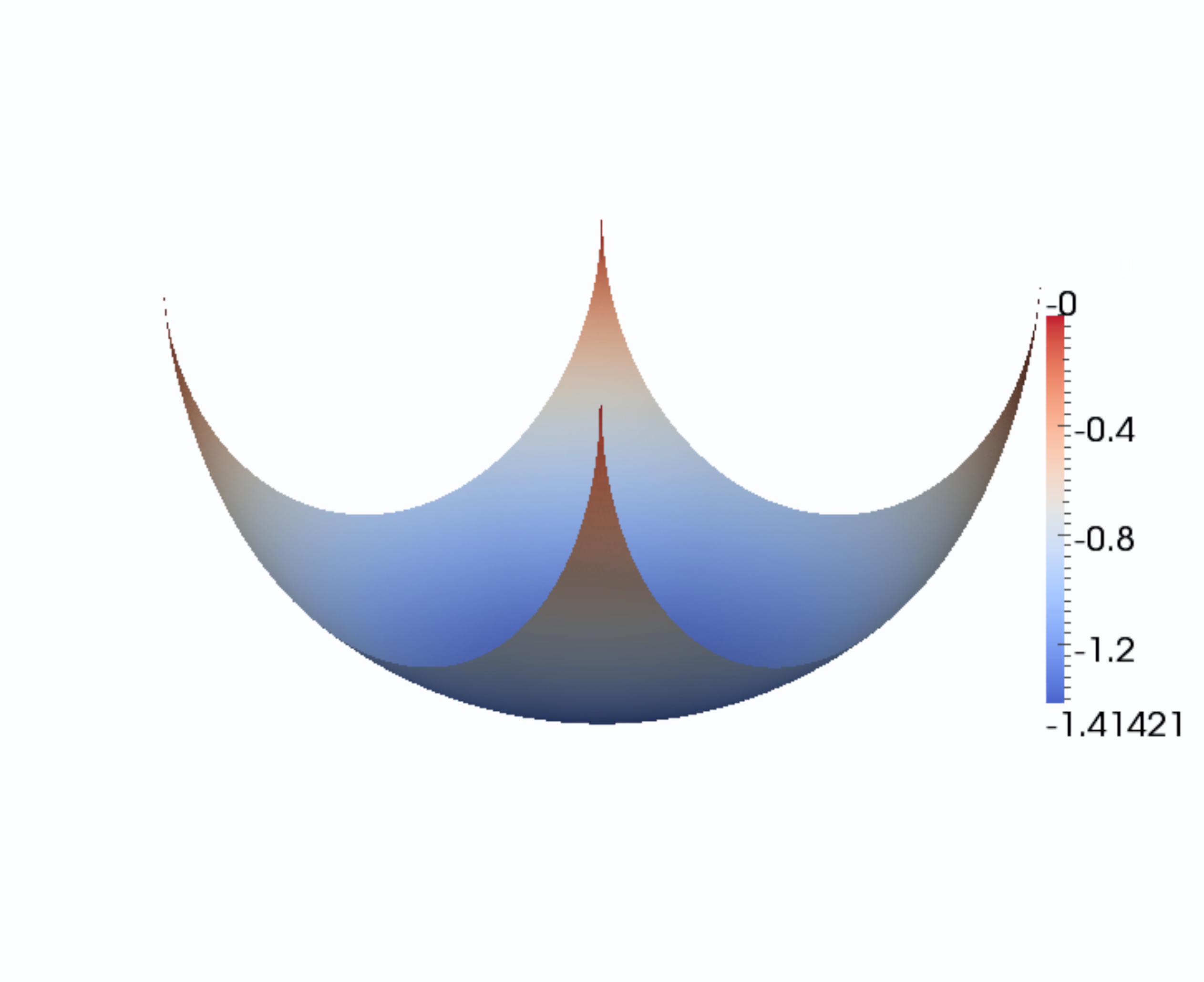}
      } \hfill 
         \subfigure[][{Log--log error plot for $\poly{2}$ Lagrange FEs.}]{
            \label{Fig:ma-ballie}
            \includegraphics[scale=\figscale,width=0.47\figwidth]{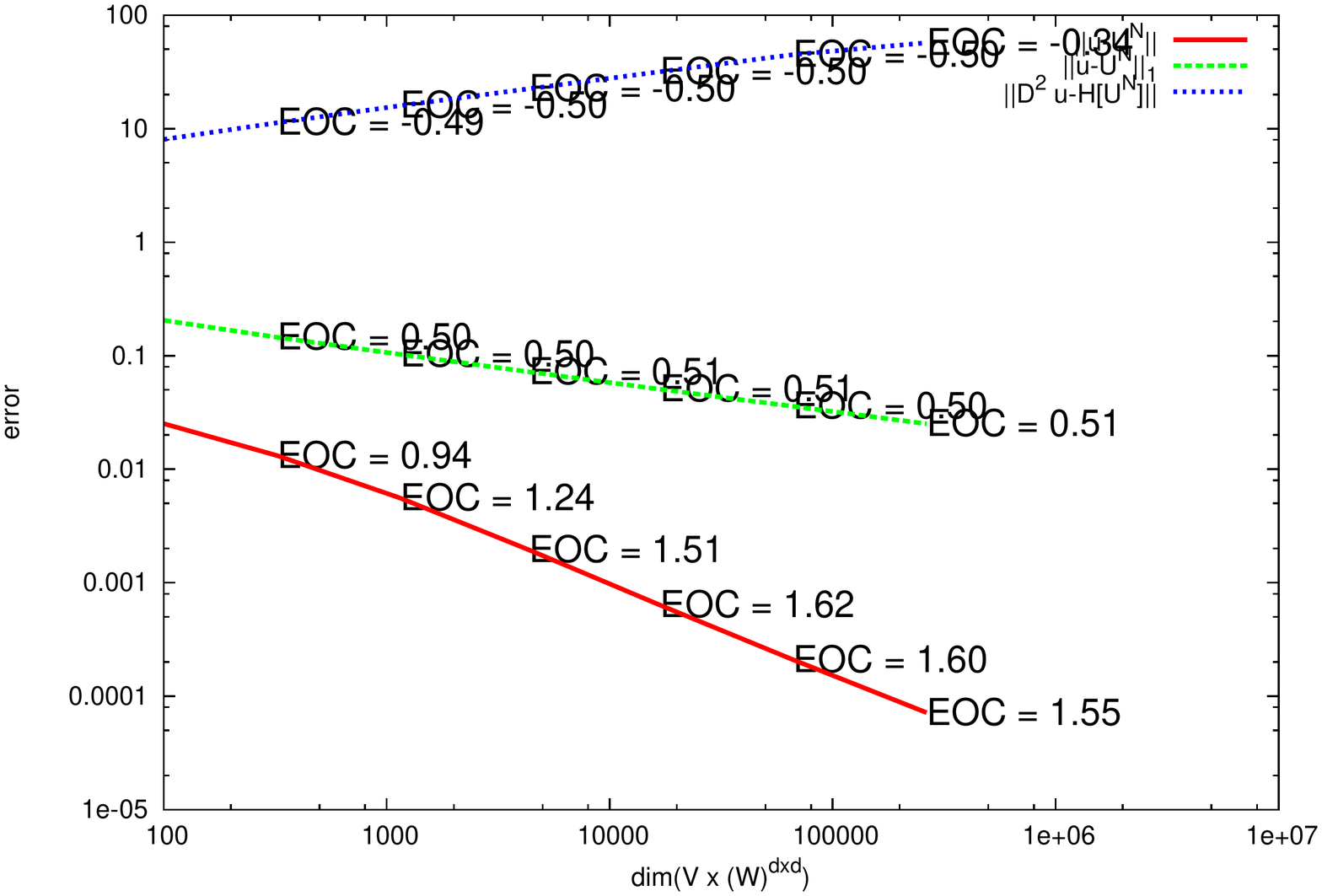}
          } 
      \end{center}
      \end{figure}

      Figures \ref{Fig:EOC-ma1}--\ref{Fig:EOC-ma2} details the various
      experiments and shows numerical convergence results for each of the
      problems studied as well as solution plots, it is worthy of note that
      each of the solutions seems to be convex, however this is not
      necessarily the case. They are all though \emph{finite element convex}
      \cite{Aguilera:2008}. In each of these cases the Dirichlet boundary
      values are not zero. The implementation of nontrivial boundary
      conditions is described in \cite[\S 3.6]{LakkisPryer:2011} or in more
      detail in \cite[\S 4.4]{Pryer:2010}.
      
      \begin{Rem}[choosing the ``right'' initial guess]
        As with any Newton method we require a starting guess, not just for
        $U^0$ but also of $\H[U^0]$. Due to the mild nonlinearity with the
        previous example an initial guess of $U^0\equiv 0$ and $\H[U^0]
        \equiv \geomat 0$ was sufficient. The initial guess to the MAD
        problem must be more carefully sought. 
        
        Since we restrict our solution to the space of convex functions, it
        is prudent for the initial guess to also be convex. Moreover we must
        rule out constant and linear functions over $\W$, since the Hessian
        of these objects would be identically zero, destroying ellipticity
        on the initial Newton step. Hence we specify that the initial guess
        to (\ref{eq:linearised-ma}) must be strictly convex.  Rather than
        postprocessing the finite element Hessian from a initial project
        (although this is an option) to initialise the algorithm we solve a
        linear problem using the nonvariational finite element
        method. Following a trick, described in \cite{DeanGlowinski:2003},
        we chose $U^0$ to be the standard $\fes$-finite element
        approximation of $u^0$ such that
        \begin{gather}
          \label{eq:laplace-initial-guess}
          \Delta u^0 = 2\sqrt{f} \text{ in } \W
          \\
          u^0 = g \text{ on } \partial\W.
        \end{gather}
      %%%%%%%%%%%%%%%%%%%%%%%%%%%%%%%%%%%%%%%%%%%%%%%%%%%%%%%%%%%%%%%%%%%%%%%%
      %%   As already noted in \cite[Thm 3.5]{LakkisPryer:2011} for this
      %%   problem (\ref{eq:laplace-initial-guess}) the nonvariational finite
      %%   element solution will coincide with its variational sibling, albeit
      %%   with the added advantage of providing us with the finite element
      %%   Hessian of the approximation, this gives us sufficient information
      %%   to begin the algorithm.
      \end{Rem}
      
      \begin{Rem}[degree of the FE space]
        \label{rem:degree-of-fe-space}
        In the previous example the lowest order convergent scheme was found
        by taking $\fes$ to be the space of piecewise linear functions
        ($p=1$). For the MAD problem we require a higher approximation
        power, hence we take $\fes$ to be the space of piecewise
        quadratic functions, \ie $p=2$. 
      
        Although the choice of $p=1$ gives a stable scheme, convergence is
        not achieved. This can be characterised by \cite[Thm
          3.6]{Aguilera:2008} that roughly says you require more
        approximation power than what piecewise linear functions provide to
        be able to approximate all convex functions. Compare with Figure \ref{Fig:EOC-ma4}.
      \end{Rem}
      
      \subsection{Nonclassical solutions}
      \label{sec:nonclassical-monge}
\begin{figure}[h]
        \caption{\label{Fig:EOC-alpha-0-55} Numerical results for the MAD problem
          on the square $S=[-1,1]^2$. Choosing $f$ and $g$ appropriately
          such that the solution is $u(\geovec x) = \norm{\geovec x}^{2\alpha}$, with $\alpha=0.55$. Note the function is singular at the origin. We plot the finite element solution together with a log--log error plot for various error functionals as in Figure
          \ref{Fig:nonlin-regular-abs-lap}.}
        \begin{center}
          \subfigure[][{The FE approximation to the function $u(\geovec x) =
              \norm{\geovec x}^{2\alpha}$, with $\alpha=0.55$.}]{
      %      \label{Fig:ma-radial}
        \includegraphics[scale=\figscale,width=0.47\figwidth]{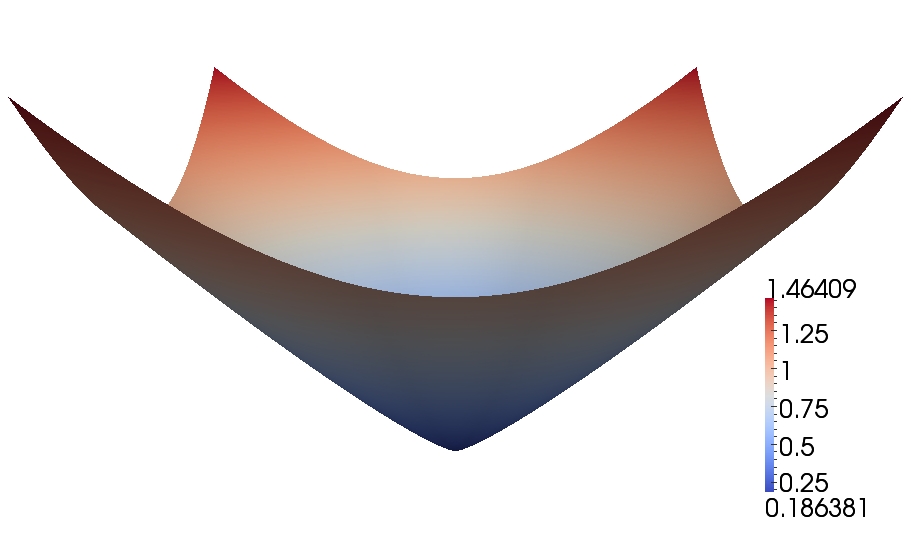}
      } \hfill 
         \subfigure[][{Log--log error plot for $\poly{2}$ Lagrange FEs.}]{
            \includegraphics[scale=\figscale,width=0.47\figwidth]{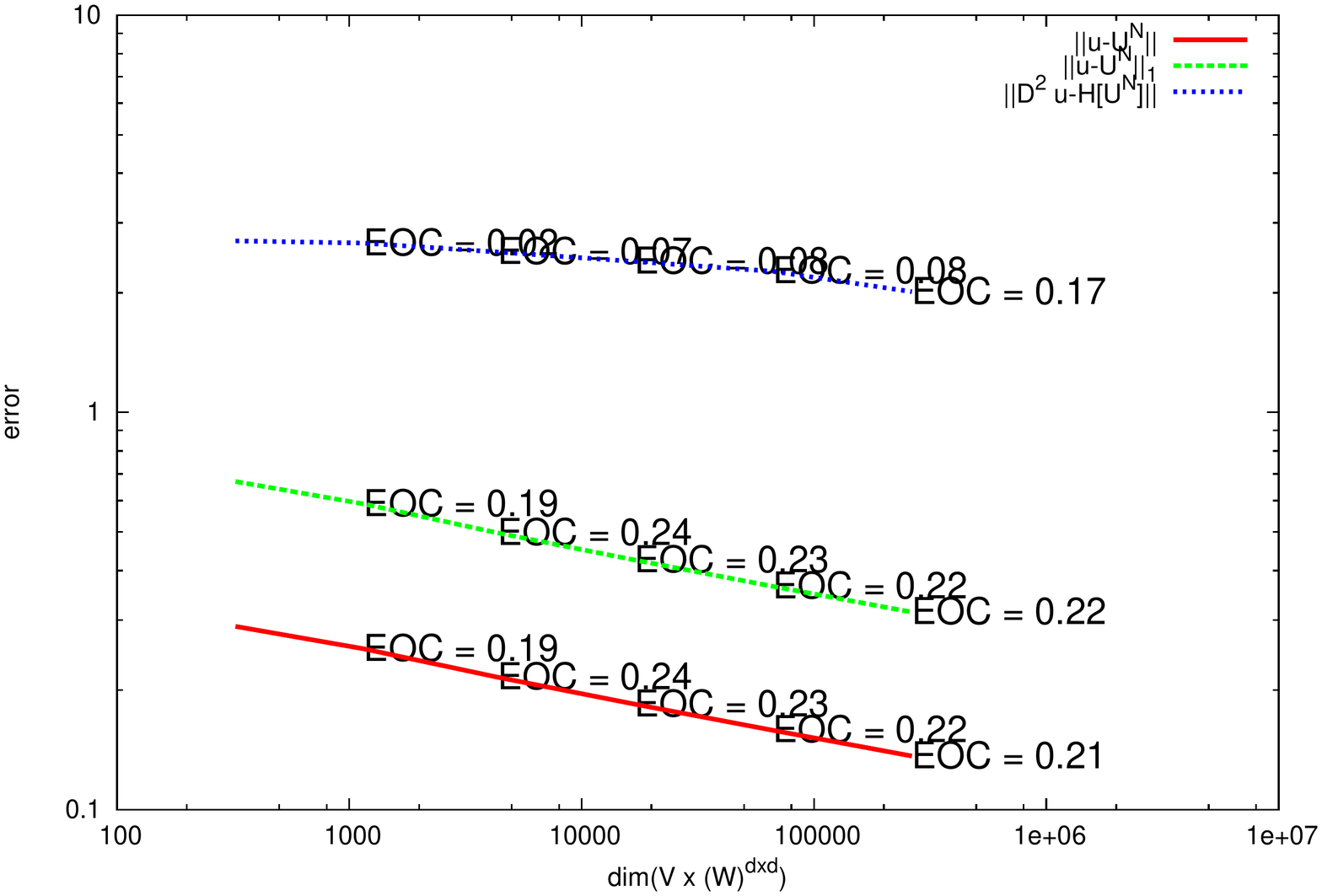}
          } 
      \end{center}
      \end{figure}

\begin{figure}[h]
        \caption{\label{Fig:EOC-alpha-0-6} Numerical results for the MAD problem
          on the square $S=[-1,1]^2$. Choosing $f$ and $g$ appropriately
          such that the solution is $u(\geovec x) = \norm{\geovec x}^{2\alpha}$, with $\alpha=0.6$. Note the function is singular at the origin. We plot the finite element solution together with a log--log error plot for various error functionals as in Figure
          \ref{Fig:nonlin-regular-abs-lap}.}
        \begin{center}
          \subfigure[][{The FE approximation to the function $u(\geovec x) =
              \norm{\geovec x}^{2\alpha}$, with $\alpha=0.6$.}]{
      %      \label{Fig:ma-radial}
        \includegraphics[scale=\figscale,width=0.47\figwidth]{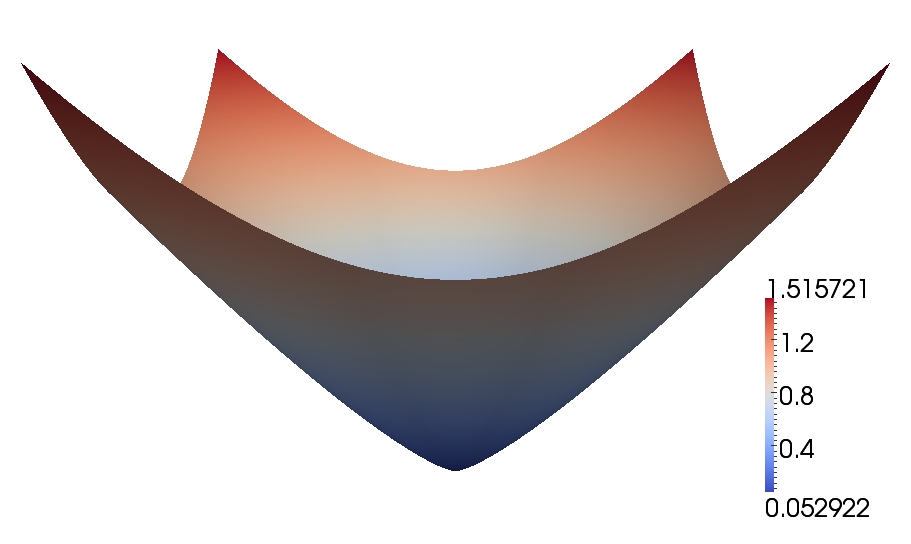}
      } \hfill 
         \subfigure[][{Log--log error plot for $\poly{2}$ Lagrange FEs.}]{
            \includegraphics[scale=\figscale,width=0.47\figwidth]{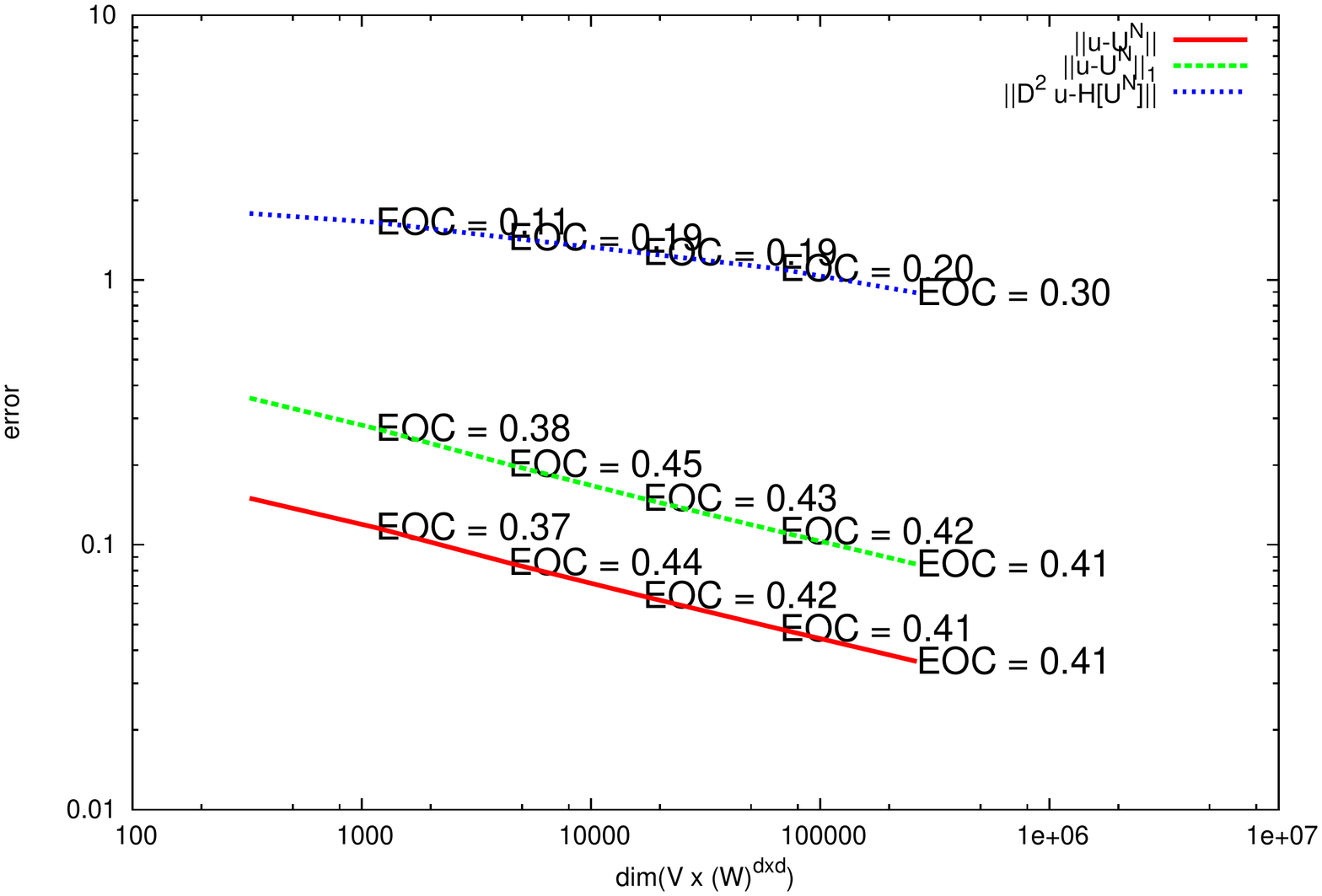}
          } 
      \end{center}
      \end{figure}

\begin{figure}[h]
        \caption{\label{Fig:EOC-alpha-0-7} Numerical results for the MAD problem
          on the square $S=[-1,1]^2$. Choosing $f$ and $g$ appropriately
          such that the solution is $u(\geovec x) = \norm{\geovec x}^{2\alpha}$, with $\alpha=0.7$. Note the function is singular at the origin. We plot the finite element solution together with a log--log error plot for various error functionals as in Figure
          \ref{Fig:nonlin-regular-abs-lap}.}
        \begin{center}
          \subfigure[][{The FE approximation to the function $u(\geovec x) =
              \norm{\geovec x}^{2\alpha}$, with $\alpha=0.7$.}]{
      %      \label{Fig:ma-radial}
        \includegraphics[scale=\figscale,width=0.47\figwidth]{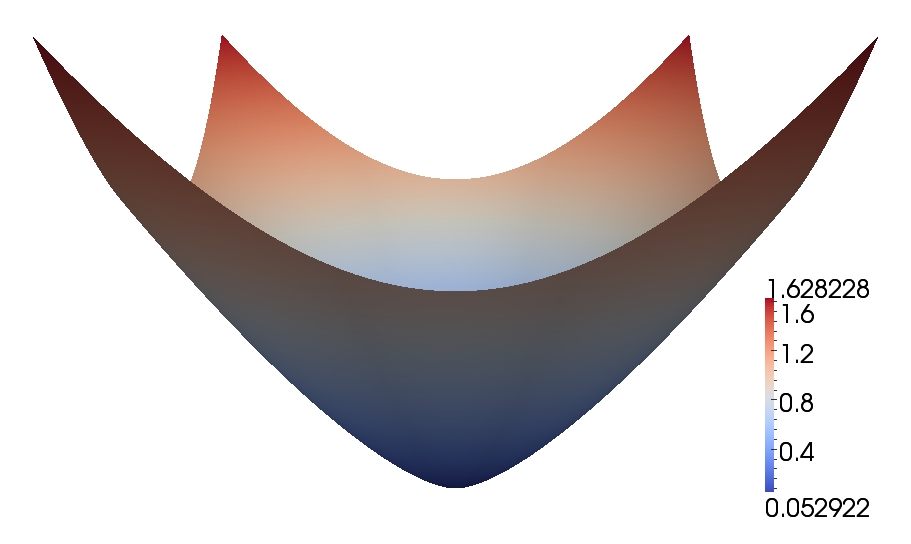}
      } \hfill 
         \subfigure[][{Log--log error plot for $\poly{2}$ Lagrange FEs.}]{
            \includegraphics[scale=\figscale,width=0.47\figwidth]{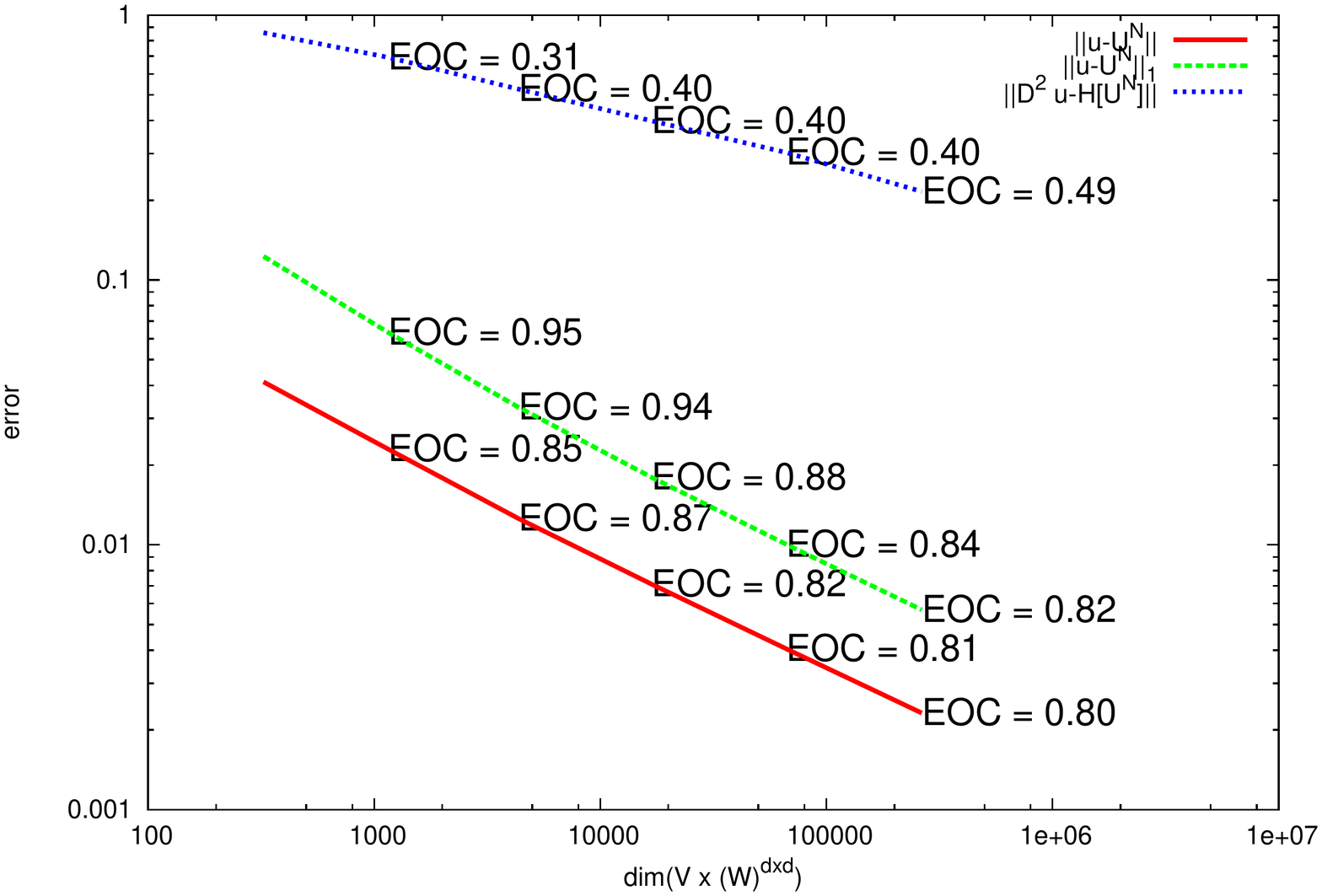}
          } 
      \end{center}
      \end{figure}

\begin{figure}[h]
        \caption[]{Numerical results for a solution to the \MAD equation with
          $f$ and $g$ appropriately
          such that the solution is $u(\geovec x) = \norm{\geovec x}^{2\alpha}$, with $\alpha=0.55$.
          We choose $p=2$, and use an adaptive
          scheme based on Z--Z gradient recovery. The mesh is refined
          correctly about the origin. Note that when $\dim{\fes} = 20,420$ the adaptive solution achieves $\Norm{u - U^M} \approx 0.0078$, the uniform solution given in Figure \ref{Fig:EOC-alpha-0-55} satisfies $\Norm{u-U^M} \approx 0.18$ using the same number of degrees of freedom. Using the adaptive strategy both $\Norm{u-U^M}$ and $\norm{u - U^M}_1$ converge like $\Oh(N^{-1})$. {\label{Fig:monge-adaptive} }}
        \begin{center}
          \subfigure[][{
              Adaptive mesh
          }]{
            \includegraphics[scale=\figscale,width=0.47\figwidth]{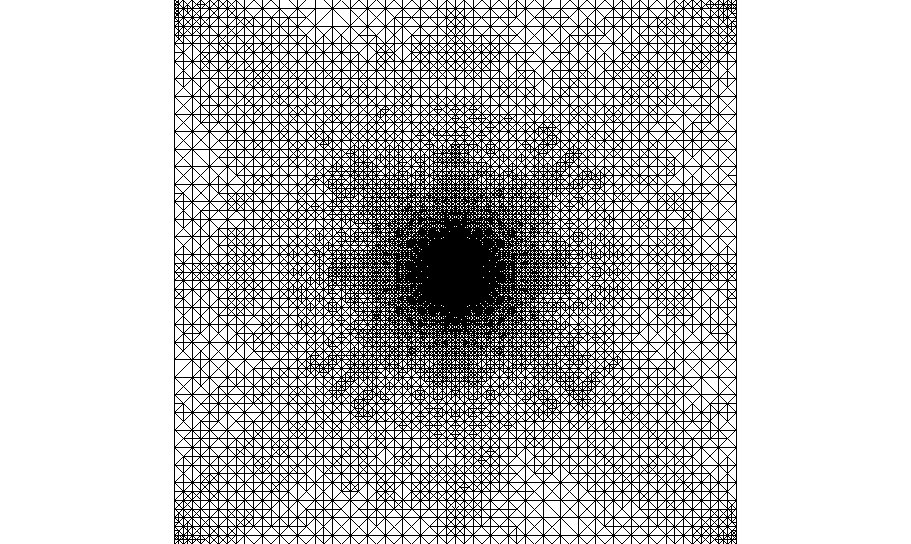}
          }
          \hfill
          \subfigure[][{Log--log error plot for $\poly{2}$ Lagrange FEs.}]{
            \includegraphics[scale=\figscale,width=0.47\figwidth]{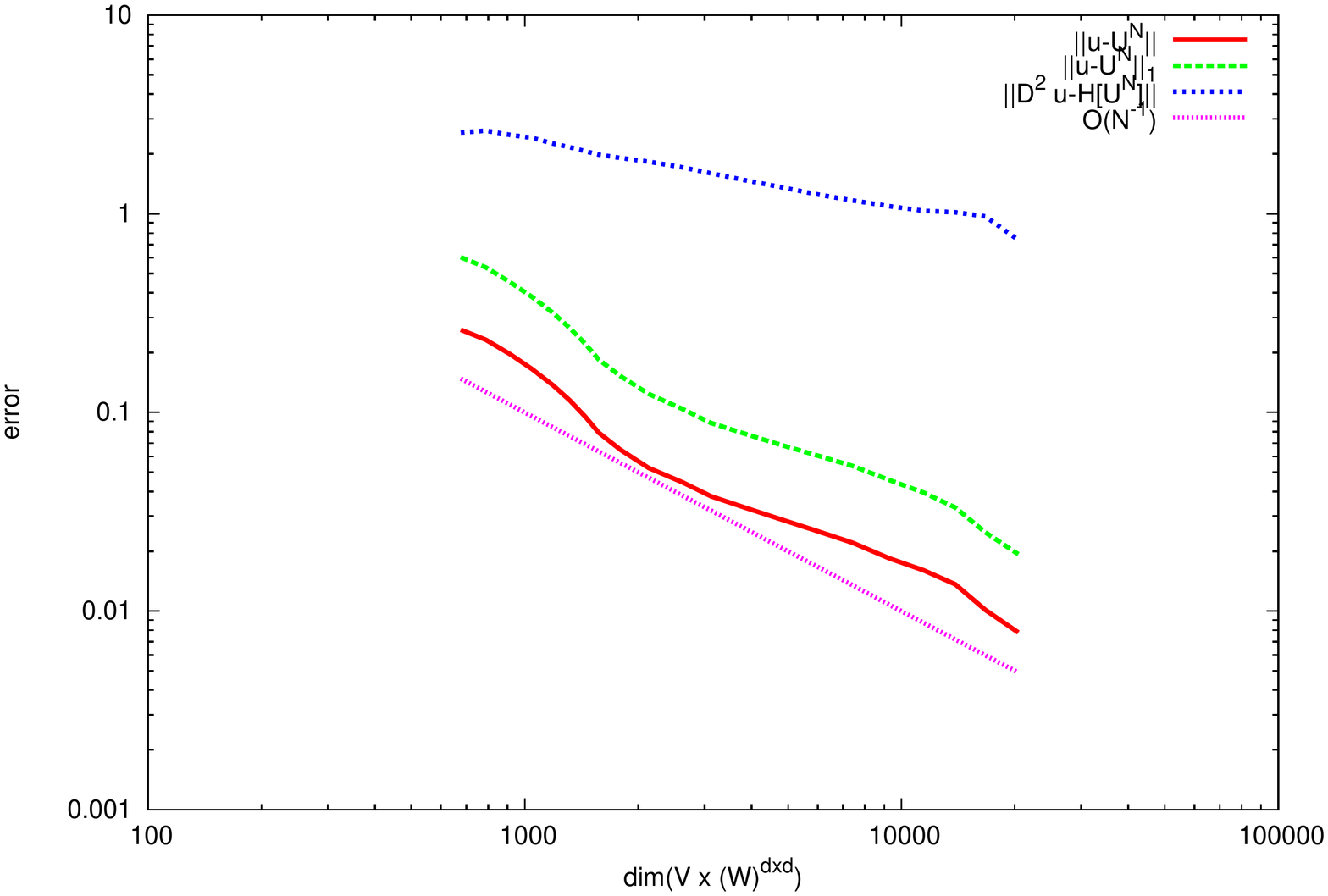}
          }
        \end{center}
      \end{figure}

      The numerical examples given in Figures \ref{Fig:EOC-ma1}--\ref{Fig:EOC-ma2} both describe the numerical approximation of classical solutions to the MAD problem. In the case of Figure \ref{Fig:EOC-ma1} $u\in\cont{\infty}(\closure\W)$ whereas in Figure \ref{Fig:EOC-ma2} $u\in\cont{\infty}(\W)\cap\cont{0}(\closure{\W})$. We now take a moment to study less regular solutions, \ie viscocity solutions which are not classical. In this test we the solution
      \begin{equation}
        u(\geovec x) = \norm{\geovec x}^{2\alpha}
      \end{equation}
      for $\alpha\in (1/2, 3/4)$. The solution $u(\geovec x) \notin \sobh{2}(\W)$. In Figures \ref{Fig:EOC-alpha-0-55}--\ref{Fig:EOC-alpha-0-7} we vary the value of $\alpha$ and study the convergence properties of the method. We note that the method fails to find a solution for $\alpha \leq 1/2$. Finally in Figure \ref{Fig:monge-adaptive} we conduct an adaptive experiment based on a gradient recovery aposteriori estimator. The recovery estimator we make use of is the Zienkiewicz--Zhu patch recovery technique see
      \cite{publication1}, \cite[\S 2.4]{Pryer:2010} or \cite[\S
        4]{Ainsworth:2000} for further details.
      
      %% \begin{figure}[h]
      %%   \caption{\label{Fig:EOC-alpha-0-8} Numerical results for the MAD problem
      %%     on the square $S=[-1,1]^2$. Choosing $f$ and $g$ appropriately
      %%     such that the solution is $u(\geovec x) = \norm{\geovec x}^{2\alpha}$, with $\alpha=0.8$. Note the function is singular at the origin. We plot the finite element solution together with a log--log error plot for various error functionals as in Figure
      %%     \ref{Fig:nonlin-regular-abs-lap}.}
      %%   \begin{center}
      %%     \subfigure[][{The FE approximation to the function $u(\geovec x) =
      %%         \norm{\geovec x}^{2\alpha}$, with $\alpha=0.8$.}]{
      %% %      \label{Fig:ma-radial}
      %%   \includegraphics[scale=\figscale,width=0.47\figwidth]{fullynonlinear/Figures/monge-singular-alpha-0-8.jpg}
      %% } \hfill 
      %%    \subfigure[][{Log--log error plot for $\poly{2}$ Lagrange FEs.}]{
      %%       \includegraphics[scale=\figscale,width=0.47\figwidth]{fullynonlinear/Figures/monge-singular-alpha-0-8-convergence}
      %%     } 
      %% \end{center}
      %% \end{figure}
      
      \clearpage
%%%%%%%%%%%%%%%%%%%%%%%%%%%%%%%%%%%%%%%%%%%%%%%%%%%%%%%%%%%%%%%%%%%%%%%%
      %%%%%%%%%%%%%%%%%%%%%%%%%%%%%%%%%%%%%%%%%%%%%%%%%%%%%%%%%%%%%%%%%%%%%%%%
      %%% mode:latex ***
      %%% tex-main-file: "../nonlinear.tex"  ***
      %%%%%%%%%%%%%%%%%%%%%%%%%%%%%%%%%%%%%%%%%%%%%%%%%%%%%%%%%%%%%%%%%%%%%%%%
      \section{Pucci's equation}
      \label{sec:Pucci}
      In this section we look to discretise the nonlinear problem, in this
      case Pucci's equation as a system of nonlinear equations.  Pucci's
      equation arises as a linear combination of Pucci's extremal
      operators. It can nevertheless be written in an algebraically
      accessible form, without the need to compute the eigenvalues.
      %%%%%%%%%%%%%%%%%%%%%%%%%%%%%%%%%%%%%%%%%%%%%%%%%%%%%%%%%%%%%%%%%%%%%%%%
      \begin{Defn}[Pucci's extremal operators \cite{CaffarelliCabre:1995}]
        Let $\geomat N \in \symm$ and $\eig{\geomat N}$ be it's spectrum,
        then the extremal operators are
        \begin{gather}
          \cM(\geomat N) := \sum_{\lambda_i \in \sigma(\geomat N)} \alpha_i \lambda_i = 0,
        \end{gather}
        with $\alpha_i\in\reals$. The maximal (minimal) operator, commonly
        denoted $\cM^+$ ($\cM^-$), has coefficients that satisfy
        \begin{equation}
          0 < \alpha_1 \leq \dots \leq \alpha_n 
          \qquad
          \qp{ \alpha_1 \geq \dots \geq \alpha_n > 0}
        \end{equation}
        respectively.
      \end{Defn}
      %%%%%%%%%%%%%%%%%%%%%%%%%%%%%%%%%%%%%%%%%%%%%%%%%%%%%%%%%%%%%%%%%%%%%%%%
      \subsection{The planar case and uniform ellipticity}
      In the case $d=2$ the normalised Pucci's equation reduces to finding $u$ such that
      \begin{equation}
        \label{eq:Pucci}
        \alpha \lambda_2 + \lambda_1 = 0
      \end{equation}
      where $\geomat N := \Hess u$. Note that if $\alpha = 1$
      (\ref{eq:Pucci}) reduces to the Poisson--Dirichlet problem. This can
      be easily seen when reformulating the problem as a second order PDE
      \cite{DeanGlowinski:2005}. Making use of the characteristic
      polynomial, we see
      \begin{equation}
        \begin{split}
          \lambda_i = \frac{\Delta u \pm \qp{\qp{\Delta u}^2 - 4\det{\Hess u}}^{1/2}}{2} \qquad i = 1,2.
        \end{split}
      \end{equation}
      Thus Pucci's equation can be written as
      \begin{equation}
        \label{eq:pucci-pde}
        0 = \qp{\alpha + 1}\Delta u + \qp{\alpha - 1}\qp{\qp{\Delta u}^2 - 4\det{\Hess u}}^{1/2},
      \end{equation}
      which is a nonlinear combination of \MA and Poisson problems.  However
      owing to the Laplacian terms, and unlike the \MAD problem, Pucci's
      equation is (unconditionally) uniformly elliptic for
      %%\begin{Pro}
      %%  Given a matrix 
        \begin{equation}
          \qp{\tr{\geomat X}}^2 - 4\det{\geomat X} 
          \geq
          0
          \Foreach
          \geomat X \in \reals^{2\times 2}.
        \end{equation}
      %%\end{Pro}
      %% \begin{Proof}
      %%   The proof is elementary linear algebra. Let $\lambda_{1,2}$ denote
      %%   the eigenvalues of $\geomat X$. Writing the trace and determinant as
      %%   functions of the eigenvalues gives
      %%   \begin{equation}
      %%     \begin{split}
      %%       \qp{\tr{\geomat X}}^2 - 4\det{\geomat X} 
      %%       &=
      %%       \qp{\sum_{i=1}^2 \lambda_i}^2 - 4\prod_{i=1}^2 \lambda_i
      %%       \\
      %%       &=
      %%       \qp{\lambda_1 - \lambda_2}^2 \geq 0,
      %%     \end{split}
      %%   \end{equation}
      %%   as required.
      %% \end{Proof}

      The discrete problem we use is a direct approximation of
      (\ref{eq:pucci-pde}), we seek $\qp{U, \H[U]}$ such that 
      \begin{gather}
          \int_\W \bigg(\qp{\alpha + 1}\tr{\H[U]} + \qp{\alpha - 1}\qp{\qp{\tr{\H[U]}}^2 - 4\det{\H[U]}}^{1/2}\bigg)\Phi = 0
          \\
          \ltwop{\H[U]}{\Psi} 
          = 
          -
          \int_\W{\nabla U}\otimes{\nabla \Psi} 
          +
          \int_{\partial \W} {\nabla U}\otimes{\geovec{n} \ \Psi}
          \Foreach (\Phi,\Psi)\in\feszero \times \fes.
      \end{gather}
      
      The result is a nonlinear system of equations which was solved using a
      algebraic Newton method.
      
      \subsection{Numerical experiments}
      
\begin{figure}
        \caption[]{Numerical results for a classical solution to Pucci's
          equation (\ref{eq:Pucci-solution}). As with the case of the MAD
          problem we choose $p=2$. We use a Newton method to solve the
          algebraic system until the residual of the problem (see
          \cite[c.f.]{Kelley:1995}) is less than $10^{-10}$ (which is
          overkill to minimise Newton error effects). We plot log--log error
          plots with experimental orders of convergence, for various norms
          and values of $\alpha$.{\label{Fig:pucci-converence} }}
        \begin{center}
          \subfigure[][{
          $\alpha = 2$
          }]{
            \includegraphics[scale=\figscale,width=0.47\figwidth]{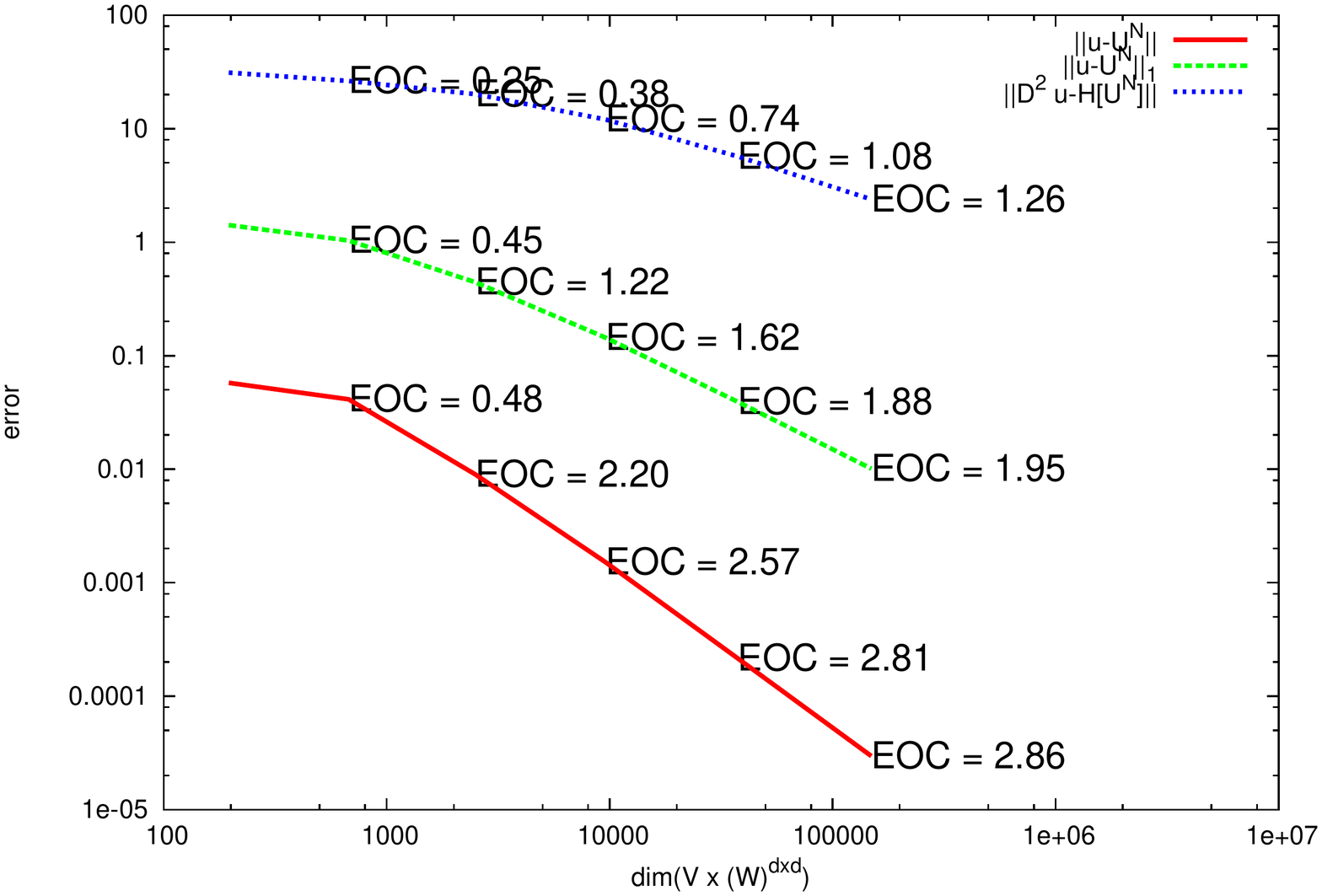}
          }
          \hfill
          \subfigure[][{$\alpha =3$}]{
            \includegraphics[scale=\figscale,width=0.47\figwidth]{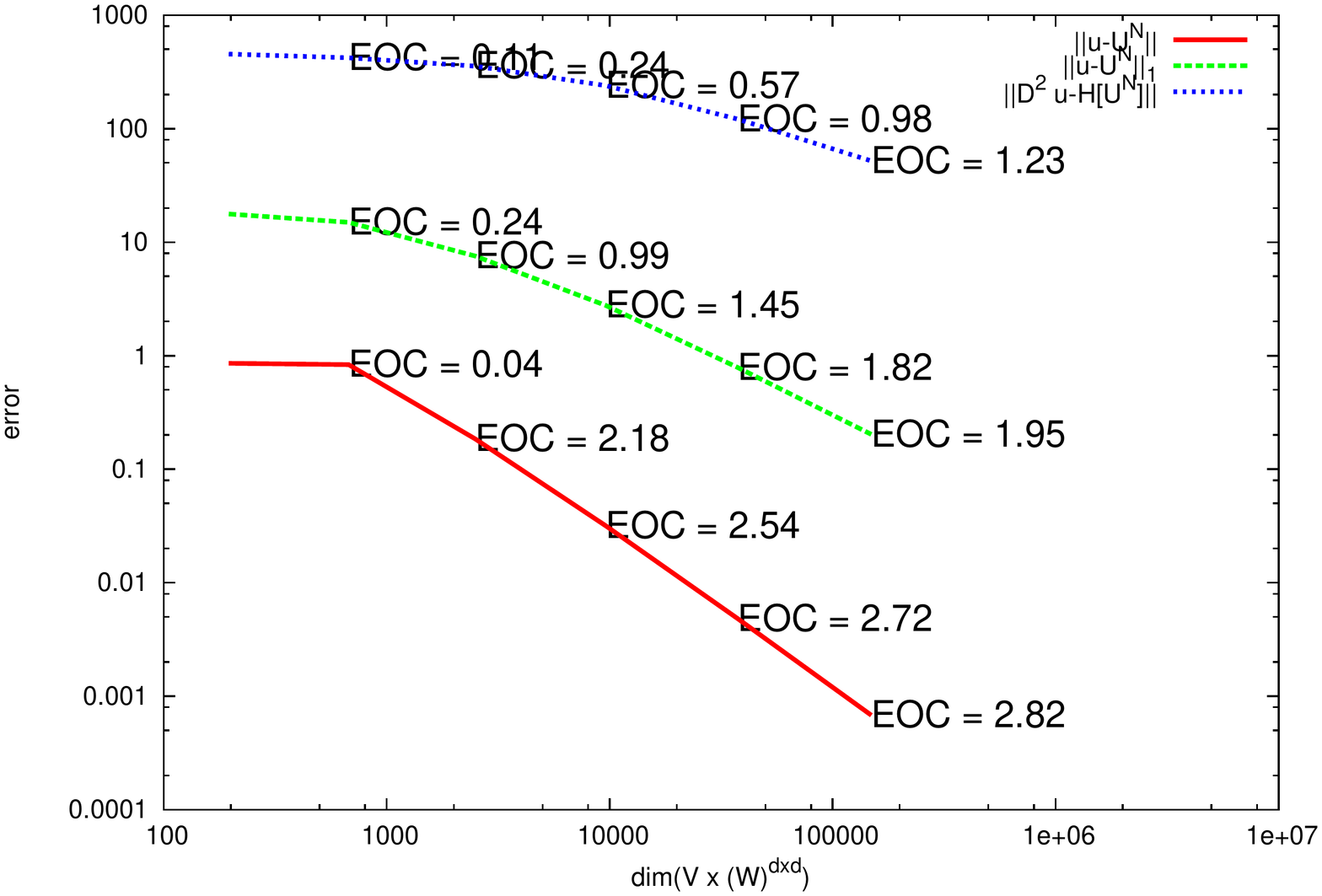}
          }
          \hfill
          \subfigure[][{
              $\alpha = 4$
          }]{
            \includegraphics[scale=\figscale,width=0.47\figwidth]{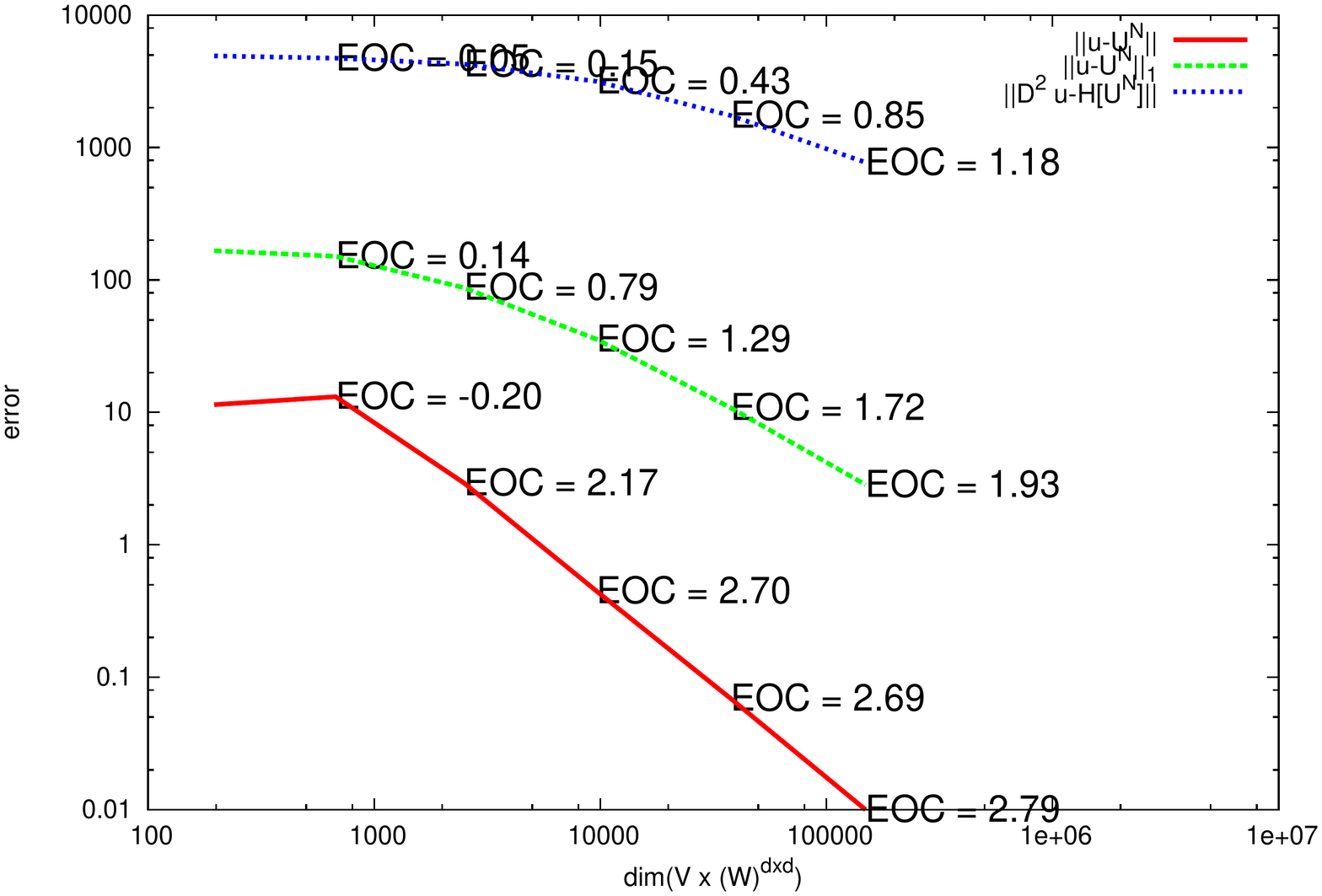}
          }
          \hfill
          \subfigure[][{$\alpha = 5$}]{
            \includegraphics[scale=\figscale,width=0.47\figwidth]{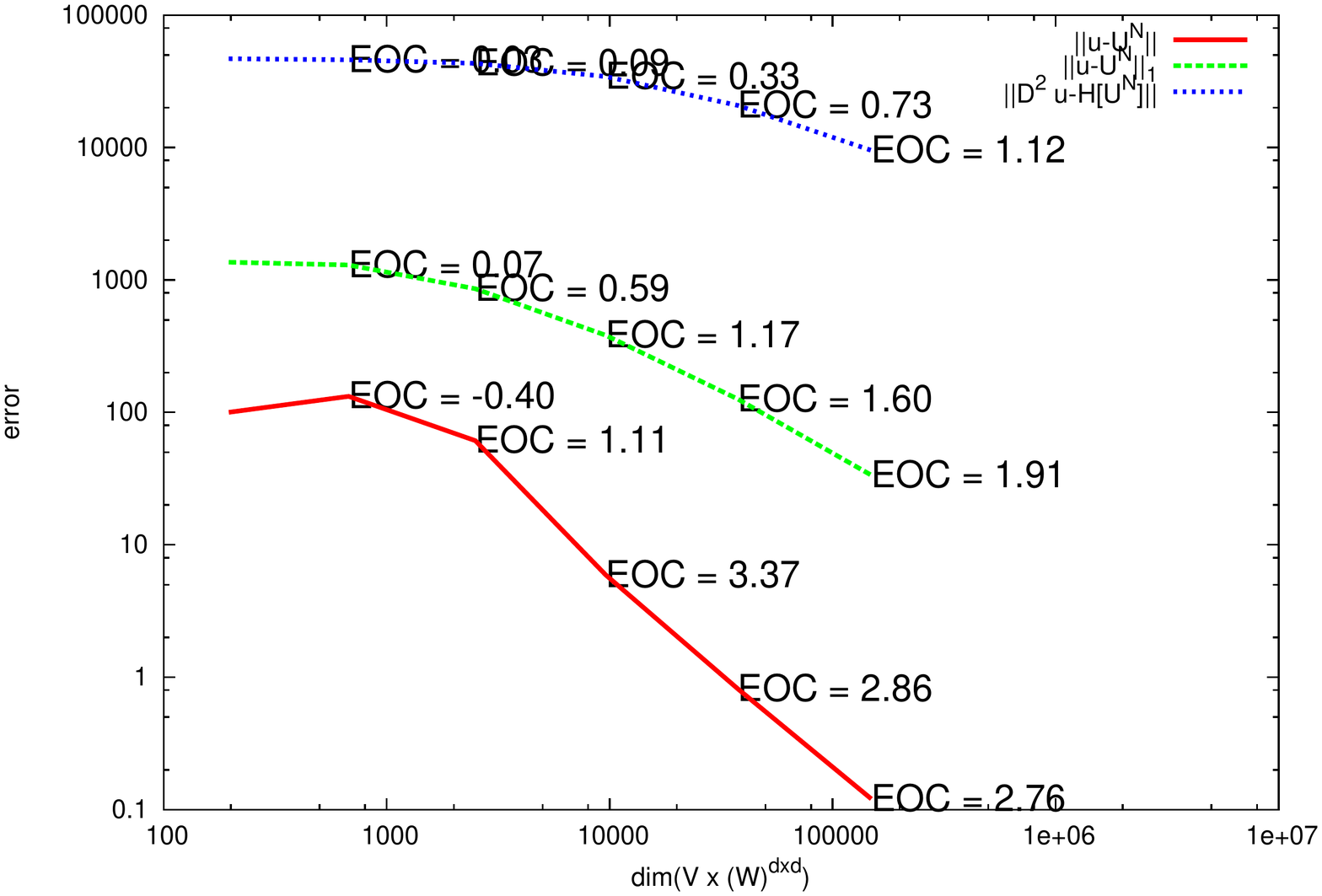}
          }
        \end{center}
      \end{figure}

      We conduct numerical experiments to be compared with those of
      \cite{DeanGlowinski:2005}. The first problem we consider is a
      classical solution of Pucci's equation (\ref{eq:Pucci}). Let $\geovec
      x = \Transpose{\qp{x, y}}$, then the function
      \begin{equation}
        \label{eq:Pucci-solution}
        u(\geovec x) = -\qp{\qp{\qp{x + 1}^2 + \qp{y+1}^2}^{\qp{1-\alpha}/2}}
      \end{equation}
      solves Pucci's equation almost everywhere away from $(x, y) = (-1,-1)$
      with $g := u|_{\partial\W}$. Let $\T{}$ be an irregular triangulation
      of $\W = [-0.95, 1]^2$. In Figure \ref{Fig:pucci-converence} we detail
      a numerical experiement considering the case $\alpha \in [2,5]$.
      
\begin{figure}[h]
        \caption[]{Numerical results for a solution to Pucci's equation with
          a piecewise defined boundary condition
          (\ref{eq:pucci-pw-boundary}). We choose $p=2$, use a Newton method
          to solve the algebraic system until the residual of the problem is
          less than $10^{-10}$. We plot the solution for various values of
          $\alpha$ as well as a cross section through the coordinate
          axis. Notice that the solution becomes extremely badly behaved as
          $\alpha$ increases. {\label{Fig:pucci-piecewise} }}
        \begin{center}
          \subfigure[][{
          $\alpha = 2$
          }]{
            \includegraphics[scale=\figscale,width=0.47\figwidth]{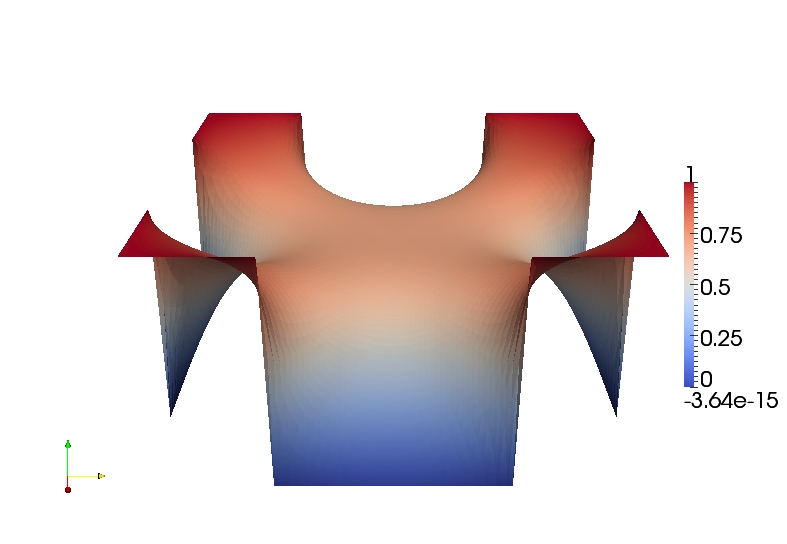}
          }
          \hfill
          \subfigure[][{$\alpha =3$}]{
            \includegraphics[scale=\figscale,width=0.47\figwidth]{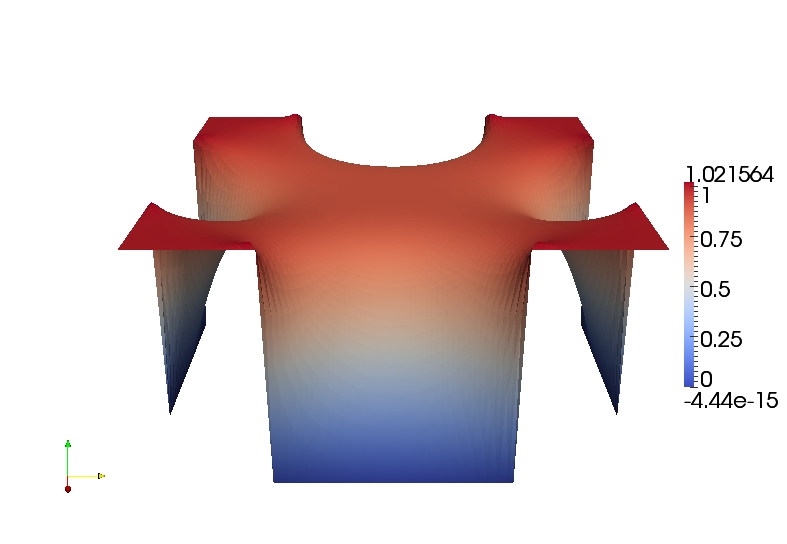}
          }
          \hfill
          \subfigure[][{
              $\alpha = 4$
          }]{
            \includegraphics[scale=\figscale,width=0.47\figwidth]{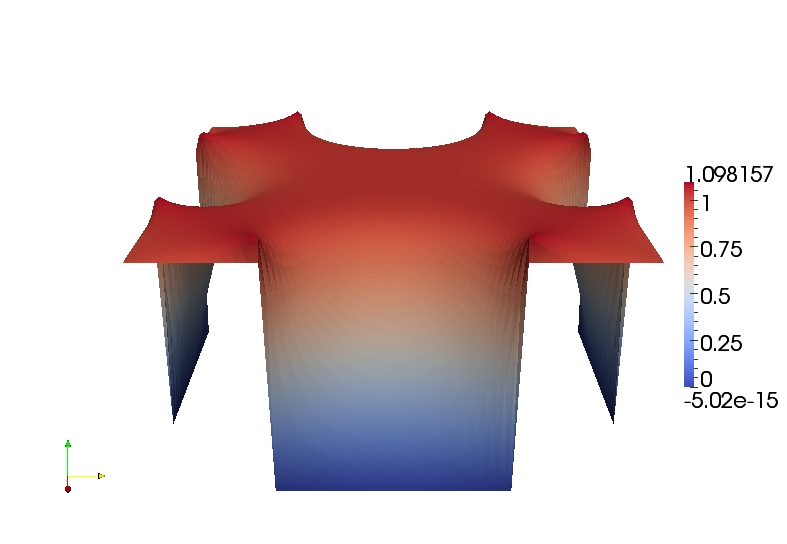}
          }
          \hfill
          \subfigure[][{$\alpha = 5$}]{
            \includegraphics[scale=\figscale,width=0.47\figwidth]{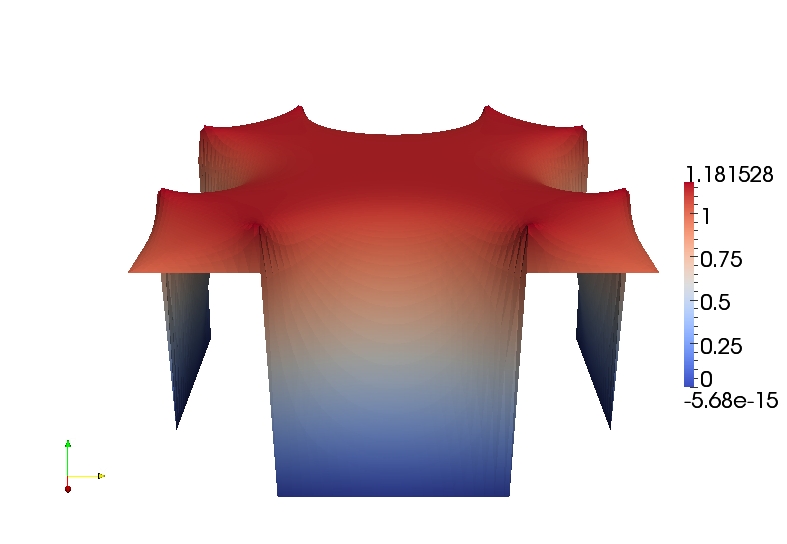}
          }
          \hfill
          \subfigure[][{Cross section about $x = 0$ for each value of $\alpha$.}]{
            \includegraphics[scale=\figscale,width=0.47\figwidth]{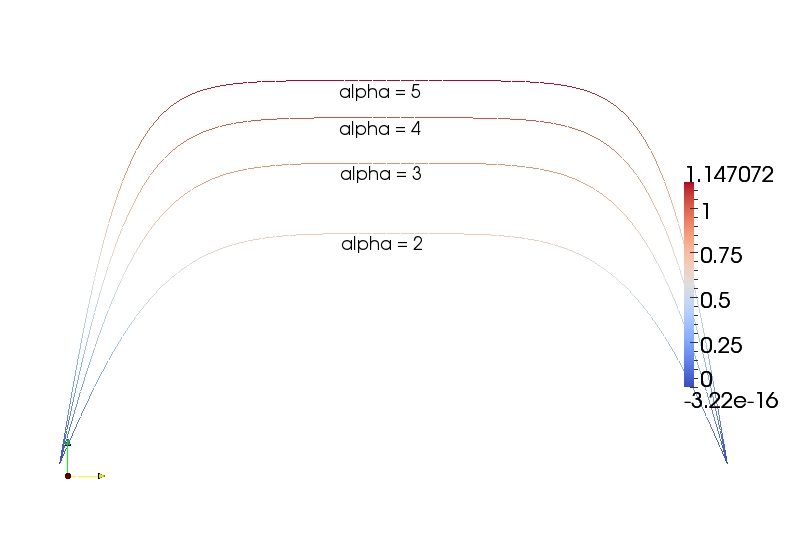}
          }
        \end{center}
      \end{figure}

      We also conduct a numerical experiment to be compared with
      \cite{Oberman:2008}. In this problem we consider a solution of Pucci's
      equation with a piecewise defined boundary. Let $\W = [-1,1]^2$ and
      the boundary data be given as
      \begin{equation}
        \label{eq:pucci-pw-boundary}
        g(\geovec x):= 
        \begin{cases}
          1 \qquad \text{ when } \abs{x} \geq \frac{1}{2} \AND \abs{y} \geq \frac{1}{2}
          \\
          0 \qquad \text{ otherwise.}
        \end{cases}
      \end{equation}
      Figure \ref{Fig:pucci-piecewise} details the numerical experiment on
      this problem with various values of $\alpha$.
      
\begin{figure}[h]
        \caption[]{Numerical results for a solution to Pucci's equation with
          a piecewise defined boundary condition
          (\ref{eq:pucci-pw-boundary}). We choose $p=2$, and use an adaptive
          scheme based on Z--Z gradient recovery. The mesh is refined
          correctly about the jumps on the
          boundary. {\label{Fig:pucci-adaptive} }}
        \begin{center}
          \subfigure[][{
          Finite element solution for $\alpha = 3$
          }]{
            \includegraphics[scale=\figscale,width=0.47\figwidth]{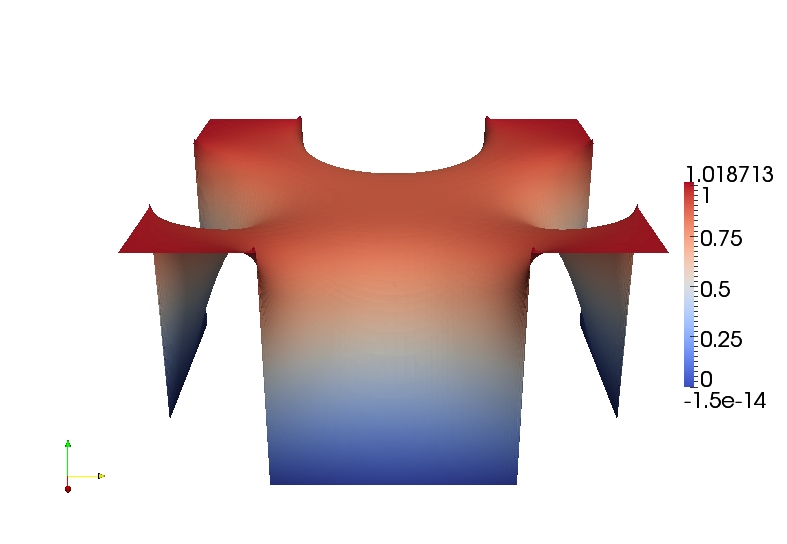}
          }
          \hfill
          \subfigure[][{Adaptive mesh}]{
            \includegraphics[scale=\figscale,width=0.47\figwidth]{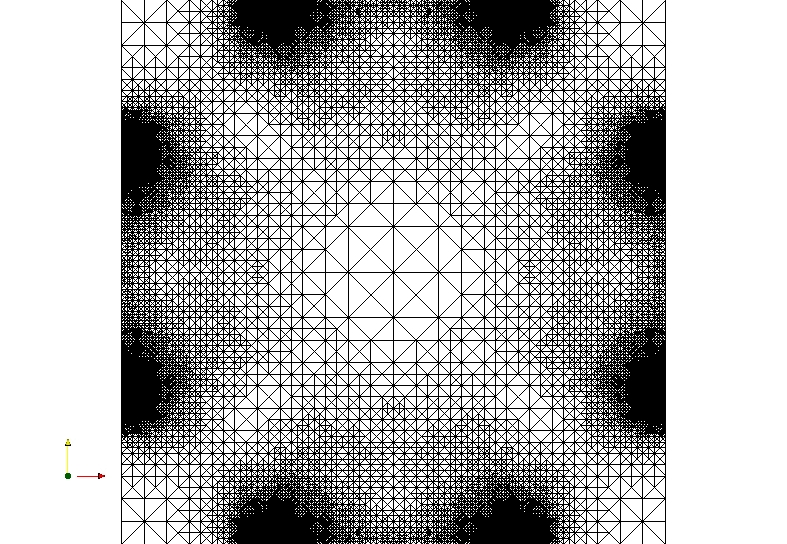}
          }
        \end{center}
      \end{figure}

      Since the solution to the Pucci's equation with piecewise boundary
      (\ref{eq:pucci-pw-boundary}) is clearly singular near the
      discontinuities we have also conducted an adaptive experiment based on
      a gradient recovery aposteriori estimator (as in \S\ref{sec:nonclassical-monge}). As can be seen from Figure
      \ref{Fig:pucci-adaptive} we regain qualitively similar results using
      far fewer degrees of freedom.
      
      %%%%%%%%%%%%%%%%%%%%%%%%%%%%%%%%%%%%%%%%%%%%%%%%%%%%%%%%%%%%%%%%%%%%%%%%
      %%%
      %%% Local Variables: ***
      %%% mode:latex ***
      %%% tex-main-file: "fullynl.tex"  ***
      %%% End: ***
%%%%%%%%%%%%%%%%%%%%%%%%%%%%%%%%%%%%%%%%%%%%%%%%%%%%%%%%%%%%%%%%%%%%%%%%
      \section{Conclusions and outlook}
      
      In this work we have proposed a novel numerical scheme for fully
      nonlinear and generic quasilinear PDEs. The scheme was based on a
      previous work for nonvariational PDEs (those given in nondivergence
      form) \cite{LakkisPryer:2011}.
      
      We have illustrated the application of the method for a simple, non
      physically motivated example, moving on to more interesting problems,
      that of the \MA equation and Pucci's equation. 
      
      For Pucci's equation we numerically showed convergence and conducted
      experiments which may be compared with previous numerical studies.
      
      We demonstrated for classical solutions to the \MA equation the method
      is robust again showing numerical convergence. For less regular
      viscocity solutions we have found that the method must be augmented
      with a penalty term in a similar light to
      \cite{BrennerGudiNeilanSung:2011}.
      
      We postulate that the method is better suited to a discontinuous
      Galerkin framework which is the subject of ongoing research.
      
      %%%%%%%%%%%%%%%%%%%%%%%%%%%%%%%%%%%%%%%%%%%%%%%%%%%%%%%%%%%%%%%%%%%%%%%%
      %%%
      %%% Local Variables: ***
      %%% mode:latex ***
      %%% tex-main-file: "fullynl.tex"  ***
      %%% End: ***
%%%%%%%%%%%%%%%%%%%%%%%%%%%%%%%%%%%%%%%%%%%%%%%%%%%%%%%%%%%%%%%%%%%%%%%%
\bibliographystyle{abbrvnat}
              
%%%%%%%%%%%%%%%%%%%%%%%%%%%%%%%%%%%%%%%%%%%%%%%%%%%%%%%%%%%%%%%%%%%%%%%%
%\printglossary  %% Print glossary
\end{document}